\documentclass[leqno]{article}
\usepackage{amsfonts,amsmath,amssymb,amsthm}

\usepackage{graphicx,epsfig}
\usepackage{checkend}
\usepackage{longtable,geometry}
\newtheorem{definition}{Definition}[section]
\newtheorem{theorem}{Theorem}
\newtheorem{lemma}[theorem]{Lemma}
\numberwithin{equation}{section}

\begin{document}

%%%%%%%%%%%%%%%%%%
\begin{center}  {\huge\textbf{Evolution of Interfaces for the Nonlinear Parabolic p-Laplacian Type Reaction-Diffusion Equations }}
%\par\medskip\bigskip \huge\textsc{Project Description}
\par \medskip\bigskip\end{center}
\begin{center} {\Large\textsc{Ugur G. Abdulla and Roqia Jeli}}
\par \medskip\bigskip\end{center}
\begin{center} {\large\noindent \textsc{Department of Mathematics, Florida Institute of Technology, Melbourne, Florida 32901}}
\par \medskip\bigskip\end{center}
{\bf Abstract.} We present a full classification of the short-time behaviour of the interfaces and local solutions to the nonlinear parabolic $p$-Laplacian type reaction-diffusion equation of non-Newtonian elastic filtration
\[ u_t-\Big(|u_x|^{p-2}u_x\Big)_x+bu^{\beta}=0, \  p>2, \beta >0 \]
The interface may expand, shrink, or remain stationary as a result of the competition of the diffusion and reaction terms near the interface, expressed in terms of the parameters $p,\beta, sign~b$, and asymptotics of the initial function near its support. In all cases, we prove the explicit formula for the interface and the local solution with accuracy up to constant coefficients. The methods of the proof are based on nonlinear scaling laws, and a barrier technique using special comparison theorems in irregular domains with characteristic boundary curves.

%{\bf Key words:} nonlinear degenerate parabolic PDE, parabolic $p$-Laplacian, reaction-diffusion equation, interface, nonlinear scaling laws, super- and %subsolutions

%{\bf AMS subject classifications:} 35K55, 35K65
%\newpage

\section{Inrtroduction}
We consider the Cauchy problem(CP) for the nonlinear degenerate parabolic equation
\begin{equation}\label{CP1}
Lu\equiv u_t-\Big(|u_x|^{p-2}u_x\Big)_x+bu^{\beta}=0, \ x\in \mathbb{R}, 0<t<T,
\end{equation}
with
\begin{equation}\label{CP2}
u(x,0)=u_0(x),~~x\in \mathbb{R},
\end{equation}
where~~$p>2,~b\in \mathbb{R},~\beta>0,~0<T\leq +\infty,$~~and~$u_0$~is nonnegative and continuous. We assume that ~$b>0$~if~$\beta <1$,~and~$b$ is arbitrary if ~$\beta \geq 1$~(see Remark 1.1). Equation \eqref{CP1} arises in many applications, such as the filtration of non-Newtonian fluids in porous media (\cite{Barenblatt1} or nonlinear heat conduction (\cite{Barenblatt2}) in the presence of the reaction term expressing additional release ($b>0$) or absorption ($b<0$) of energy. 

The goal of this paper is to analyze the behavior of interfaces separating the regions where ~$u=0$~and where ~$u>0$. We present full classification of the  short-time evolution of interfaces and local structure of solutions near the interface. Due to invariance of \eqref{CP1} with respect to translation, without loss of generality, we will investigate the case when ~$\eta(0)=0,$~where

\[\eta(t)=\text{sup}~\{x:u(x,t)>0\}.\]
and precisely, we are interested in the short-time behavior of the interface function ~$\eta(t)$~and local solution near the interface. We shall assume that 
\begin{equation}\label{CP3}
u_0\sim C(-x)_+^{\alpha}~~\text{as}~x\rightarrow 0-~~~\text{for some}~~~C>0,~\alpha>0.\end{equation}
The direction of the movement of the interface and its asymptotics is an outcome of the competition between the diffusion and reaction terms and depends on the parameters ~$p,b,\beta, C,$~and~$\alpha$. Since the main results are local in nature, without loss of generality  we may suppose that ~$u_0$~either is bounded or satisfies some restriction on its growth rate as ~$x\rightarrow -\infty$~which is suitable for existence, uniqueness, and comparison results (see section 3). The special global case
\begin{equation}\label{CP4}
u_0(x)=C(-x)_+^{\alpha},\quad x\in \mathbb{R},
\end{equation}
will be considered when the solution to the problem \eqref{CP1}, \eqref{CP4} is of self-similar form. Our estimations are global in time in these special cases. 

Initial development of interfaces and structure of local solution near the interfaces is very well understood in the case of the reaction-diffusion equations with porous medium type diffusion term:
\begin{equation}\label{CP5}
u_t-(u^m)_{xx}+bu^\beta=0 \ x\in \mathbb{R}, 0<t<T,
\end{equation}
Full classification of the evolution of interfaces and the local behaviour of solutions near the interfaces in CP \eqref{CP5}, \eqref{CP2}, \eqref{CP3} was presented in \cite{Abdulla1} for the case of fast diffusion ($m>1$) case, and in \cite{Abdulla2} for the slow diffusion ($0<m<1$) case. The major obstacle in solving the interface development problem for nonlinear degenerate parabolic equations is a problem of non-uniform asymptotics in the sense of singular perturbations theory, namely that the dominant balance as ~$t\rightarrow 0+$~between the terms in \eqref{CP1}, \eqref{CP5} on curves which approach the boundary of the support on the initial line depending on how they do so. The general theory, including existence, boundary regularity, uniqueness and comparison theorems, for the reaction diffusion equations of type \eqref{CP5} in general non-cylindrical and non-smooth domains is developed in \cite{Abdulla3} in the one-dimensional case, and in \cite{Abdulla5, Abdulla6, Abdulla7} in the multi-dimensional case. Comparison theorems proved in \cite{Abdulla3} were essential tools in developing the rigorous proof method in \cite{Abdulla1, Abdulla2} for solving interface problem for the reaction-diffusion equation \eqref{CP5}.   The rigorous proof method developed in \cite{Abdulla1, Abdulla2} is based on a barrier technique using special comparison theorems in irregular domains with characteristic boundary curves.
In this paper we apply the method developed in \cite{Abdulla1} to solve the interface problem for the PDE \eqref{CP1}.

The structure of the paper is as follows: In section 2 we outline the main results. In section 3 we apply rescaling and prove for some preliminary estimations which are necessary for using our barrier technique. Finally in section 4 we prove the results of section 2.
To avoid technical difficulties, we give explicit values of some of lengthy constants in the appendix. 

{\bf Remark 1.1}. We are not interested in the special case $p=2$ of semilinear heat equation. This case was completed in \cite{Grundy1, Grundy2} (see also \cite{Abdulla1}). However, we will mention when our results extend to the limit case $p=2$. In general, the case $p=2$ is in some sense a singular limit. For example, if $b>0, 0<\beta<1, \alpha<\frac{p}{p-1-\beta}$, then we prove that the interface initially expands and 
\[\eta(t)\sim C_1t^{1/(p-\alpha(p-2))}~~\text{as}~t\rightarrow 0+.\] By passing to the limit as $p\downarrow 2$ formally, this yields a false result. 
%\[\eta(t)=O(t^{\frac{1}{2}})\quad \text{as}~t\rightarrow 0+,\]
In fact, from \cite{Grundy2} it follows that if ~$p=2$,~then 
\[\eta(t)\sim C_2\Big (t \log\frac{1}{t}\Big )^{\frac{1}{2}}\quad \text{as}~t\rightarrow 0+\]

\section{Description of main results.} 
\begin{figure}
    \centering
    \includegraphics[width = 0.7\textwidth]{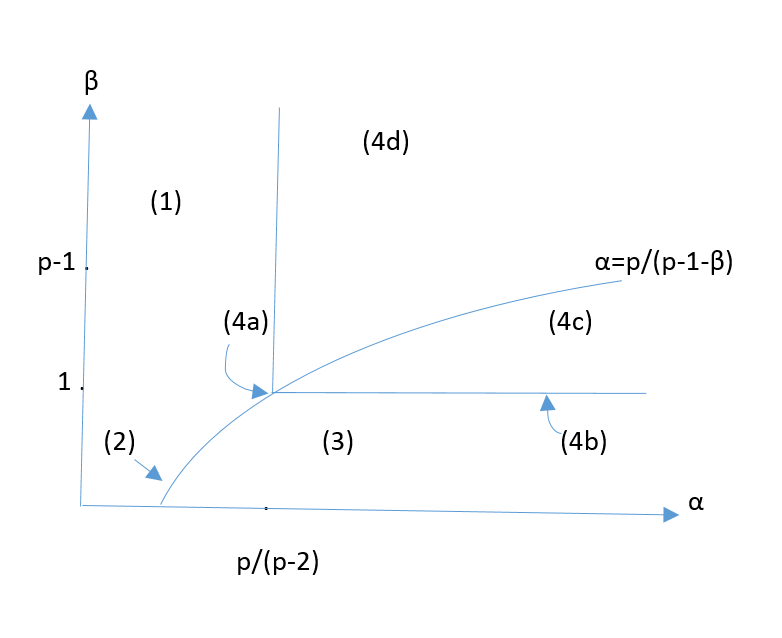}
    \caption{Classification of different cases in the ($\alpha$,$\beta$) plane for interface development in problem \eqref{CP1}-\eqref{CP4}.}
\end{figure}
In Figure 1 we present classification diagram in $(\alpha, \beta)$-plane for the initial interface development in CP \eqref{CP1}, \eqref{CP2}, \eqref{CP3} if $b>0$.\\
\begin{itemize}
\item{{\bf Region (1):}} $\alpha < p/(p-1- \min\{1,\beta\});$~Diffusion dominates and interface expands.
\item{{\bf Region (2):}} $\alpha=p/(p-1-\beta), 0<\beta <1;$~Diffusion and absorption are in balance in this borderline case. There is a critical constant $C_*$ such that interface expands for $C>C_*$, and shrinks for $C<C_*$.
\item{{\bf Region (3):}} $\alpha >p/(p-1-\beta), 0<\beta <1;$~Absorption term dominates and interface shrinks.
\item{{\bf Region (4):}} $\alpha \geq p/(p-2), \beta \geq 1;$~Interface has initial "waiting time".
\end{itemize}
To describe the asymptotic properties of the interface and local solution near the interface, we divide the results into the two different subcases:\\
{\bf(I)}    ~$b\neq 0 \  (\text{either}~ b>0, \beta>0 \ \text{or}~b<0,\beta\geq 1)$~and ~$p>2;$  \ \ \ {\bf(II)}    ~$b=0$.\\

{\bf(I)} In this case there are four different subcases, as shown in Figure 1 and itemized above. (In view of our assumptions, the case~$b<0$~relates to the part of the ($\alpha, \beta$) plane with ~$\beta\geq 1.)$\\
%\begin{itemize}

{\bf Region (1):} Let ~$\alpha < \frac{p}{p-1- \text{min}\{1,\beta\}}.$ In this case the interface initially expands and 
\begin{equation}\label{CP1'}\eta(t)\sim \xi_* t^{1/(p-\alpha(p-2))}~~\text{as}~t\rightarrow 0+,\end{equation}
where
\begin{equation}\label{CP2'}\xi_*=C^{\frac{p-2}{p-\alpha(p-2)}}\xi'_*\end{equation}
and ~$\xi'_*>0$~depends  on~$ p$~ and ~$\alpha$~only (see Lemma~\ref{lemma2}). For $\forall$ $\rho<\xi_*$ $\exists$ $f(\rho)>0$~depending on ~$C, p, $~and~$\alpha$~ such that 
\begin{equation}\label{CP3'}u(x,t)\sim f(\rho)t^{(\alpha/p-\alpha(p-2))}~~~~as~~~~t\rightarrow0+\end{equation}
along the curve ~$x=\xi_{\rho}(t)=\rho t^{1/(p-\alpha(p-2))}. $ A function $f$~is a shape function of the self-similar solution of \eqref{CP1},\eqref{CP4} with $b=0$~(see Lemma~\ref{lemma2}):
\begin{equation}\label{CP25'}u_*(x,t)=t^{\frac{\alpha}{p-\alpha(p-2)}}f(\xi),\quad ~~~\xi=xt^{-\frac{1}{p-\alpha(p-2)}},\end{equation}
In fact, $f$ is a unique solution of the following nonlinear ODE problem:
\begin{equation}\begin{cases}
\big(|f'(\xi)|^{p-2}f'(\xi)\big)'+\frac{1}{p-\alpha (p-2)}\xi f'(\xi)-\frac{\alpha}{p-\alpha (p-2)}f(\xi)=0,~-\infty <\xi <\xi_*\\
f(-\infty)\sim C(-\xi)^\alpha,f(\xi_*)=0, f(\xi)\equiv 0,~\xi \geq \xi_*
\end{cases}\label{selfsimilarODE}\end{equation}
Its dependence on $C$ is given through the following relation:
\begin{subequations}\label{24ab}
\begin{equation}\label{CP4a'} f(\rho)=C^{p/(p-\alpha(p-2))}f_0\Big(C^{(p-2)/(\alpha(p-2)-p)} \rho\Big),\end{equation}
\begin{equation}\label{CP4b'}f_0(\rho)=w(\rho,1),~~~\quad \xi'_*=\text{sup}\{\rho:f_0(\rho)>0\}>0\end{equation}\end{subequations}
where ~$w$~is a solution of \eqref{CP1}, \eqref{CP4} with ~$b=0, C=1.$~Lower and upper estimations for ~$f$~are given in \eqref{CP27'}. Moreover, 
\begin{equation}\label{CP5'}\xi'_*=A_0^{\frac{p-2}{p}}\Big[\frac{(p-1)^{p-1}(p-\alpha(p-2))}{(p-2)^{p-1}}\Big]^{\frac{1}{p}}\xi''_*,\end{equation}
where ~$A_0=w(0,1)$~and ~$\xi''_*$~is some number in ~$[\xi_1,\xi_2],$~where
\begin{equation*}
\xi_1=(p-1)^{\frac{1}{p}}\Big(\alpha(p-2)\Big)^{-\frac{1}{p}},~~~~\xi_2=1\quad~~~~~~~~~~~~~\text{if}~~~(p-1)(p-2)^{-1}\leq \alpha<p(p-2)^{-1},\end{equation*}
\begin{equation}\label{CP6'}
\xi_1=1,~~~~~~~~~~~~~~~~~\xi_2=(p-1)^{\frac{1}{p}}\Big(\alpha(p-2)\Big)^{-\frac{1}{p}},~~~~~~~~~~~~\text{if}~~~0<\alpha\leq(p-1)(p-2)^{-1}.\end{equation}
In particular, if ~$\alpha=(p-1)(p-2)^{-1}$and $p>1+(\text{min}\{1,\beta\})^{-1}$, then the explicit solution of the problem \eqref{CP1}, \eqref{CP4} with $b=0$ is given by (2.24) and we have
\begin{equation}\label{CP7'}\xi_1=\xi_2,~~~\xi'_*=(p-1)^{p-1}(p-2)^{1-p},~~~~f_0(x)=\big(\xi'_*-x\big)_+^{(p-1)/(p-2)}\end{equation}
The explicit formulae \eqref{CP1'} and \eqref{CP3'} mean that the local behavior of the interface and solution along ~$x=\xi_{\rho}(t)$~coincide with those of the problem \eqref{CP1}, \eqref{CP4} with ~$b=0$.\\

{\bf Region (2):} Let ~$b>0, 0<\beta<1, \alpha=p/(p-1-\beta)$~(here we describe the results for the case p=2 as well ). In this borderline case the direction of the movement of the interface depends on the constant $C$. The critical value is
\[C_*=\Big[\frac{|b|(p-1-\beta)^{p}}{(1+\beta)p^{p-1}(p-1)}\Big]^{\frac{1}{p-1-\beta}}\]
First, assume that $u_0$~is defined by \eqref{CP4}. If ~$\beta(p-1)=1$,~then the explicit solution to \eqref{CP1}, \eqref{CP4} is 
\begin{equation}\label{CP8'}u(x,t)=C(\zeta_*t-x)_+^{\frac{1}{1-\beta}},~~~~~\zeta_*=b(1-\beta)C^{\beta-1}((C/C_*)^{p-1-\beta}-1).\end{equation}
It has an expanding interface if ~$C>C_*$~, a shrinking interface if~$0<C<C_*$~, and is a stationary solution if~$C=C_*$.~

Let~$\beta(p-1) \neq 1.$~If~$C=C_*$~then~$u_0$~is a stationary solution to \eqref{CP1}, \eqref{CP4}. If~$C\neq C_*,$~then the solution to \eqref{CP1}, \eqref{CP4} is of the self similar form 
\begin{equation}\label{CP9'}u(x,t)=t^{1/(1-\beta)}f_1(\zeta),~~~~~~~~\zeta=xt^{-\frac{p-1-\beta}{p(1-\beta)}},\end{equation}
\begin{equation}\label{CP10'}\eta(t) =\zeta_* t^{\frac{p-1-\beta}{p(1-\beta)}},~~~~~~~~~~~~~~~~~0\leq t<+\infty.\end{equation}
If ~$C>C_*$~then the interface expands, ~$f_1(0)=A_1>0 $~(see Lemma~\ref{lemma4}), and 
\begin{equation}\label{CP11'} C_1t^{\frac{1}{1-\beta}}\Big(\zeta_1-\zeta\Big)_+^{\mu}\leq u\leq C_2t^{\frac{1}{1-\beta}}\Big(\zeta_2-\zeta\Big)_+^{\frac{p}{p-1-\beta}},~~~~0\leq x<+\infty,~~~~0<t<+\infty,\end{equation}
where
\[\mu=(p-1)(p-2)^{-1} ~~\text{if}~~ \beta(p-1)>1;~~\mu=p(p-1-\beta)^{-1}~~\text{if}~~\beta(p-1)<1\]
which implies
\begin{equation}\label{CP12'}\zeta_1\leq \zeta_*\leq \zeta_2.\end{equation}
The right-hand side of \eqref{CP11'} (respectively,\eqref{CP12'})  may be replaced by~$\bar{C_2}t^{\frac{1}{1-\beta}}(\bar{\zeta}_2-\zeta)_+^{\frac{p-1}{p-2}}$~(respectively,~$\bar{\zeta}_2$); see the appendix for the description of all the relevant constants.
Let~$\beta(p-1)\neq 1$~and~$0<C<C_*.$~Then the interface shrinks and if ~$\beta(p-1)>1$,~then \begin{equation*}\big[C^{1-\beta}(-x)_+^{\frac{p(1-\beta)}{p-1-\beta}}-b(1-\beta)t\big]_+^{\frac{1}{1-\beta}}\leq u \end{equation*}
\begin{equation}\label{CP13'}\leq \big[C^{1-\beta}(-x)_+^{\frac{p(1-\beta)}{p-1-\beta}}-b(1-\beta)(1-\Big(\frac{C}{C_*}\Big)^{p-1-\beta})t\big]_+^{\frac{1}{1-\beta}},\;~x\in \mathbb{R},\; 0\leq t<+\infty\end{equation}

which again implies \eqref{CP12'}, where $\zeta_1$(respectively,~$\zeta_2$) is replaced with
\[-C^{-\frac{p-1-\beta}{p}}\big(b(1-\beta)\big)^{\frac{p-1-\beta}{p(1-\beta)}}\]
\[\Big(\text{respectively},-C^{-\frac{p-1-\beta}{p}}(b(1-\beta)\big(1-(C/C_*)^{p-1-\beta}\big)^{\frac{p-1-\beta}{p(1-\beta)}}\Big).\]
However, if ~$\beta(p-1)<1,$~then
\begin{equation}\label{CP14'} C_*\Big(-\zeta_3t^{\frac{p-1-\beta}{p (1-\beta)}}-x\Big)_+^{\frac{p}{p-1-\beta}}\leq u \leq C_3(-\zeta_4t^{\frac{p-1-\beta}{p(1-\beta)}}-x)_+^{\frac{p}{p-1-\beta}},~0\leq t<+\infty,\end{equation}
where the left-hand side is valid for~$x\geq -\ell_0t^{\frac{p-1-\beta}{p(1-\beta)}},$ while the right-hand side is valid for ~$x\geq-\ell_1t^{\frac{p-1-\beta}{p(1-\beta)}}$. From \eqref{CP14'},\eqref{CP12'} follows if we replace ~$\zeta_1$~and~$\zeta_2$~with~$-\zeta_3$~and~$-\zeta_4$,~respectively.\\
If $\beta(p-1)\neq 1$, in general the precise value~$\zeta_*$~can be found only by solving the ODE $\mathcal{L}^0f_1=0$ (see \eqref{CP4b'''}) below) and calculating~$\zeta_*=\text{sup}~\{\zeta:f_1(\zeta)>0\}.$

Now assume that~$u_0$~satisfies \eqref{CP3} with ~$\alpha=p/(p-1-\beta.)$~Then if ~$C\neq C_*$~we have
 \begin{equation}\label{CP15'}\eta(t)\sim\zeta_*t^{\frac{p-1-\beta}{p(1-\beta)}}\quad \text{as}~~t\rightarrow 0+\end{equation}
and for $\forall \rho<\zeta_*$
\begin{equation}\label{CP16'}u(x,t)\sim f_1(\rho)t^{1/(1-\beta)}~~~~~\text{for}~~~~x=\rho t^{\frac{p-1-\beta}{p(1-\beta)}},~~~t\rightarrow 0+,\end{equation}
where the right-hand side of \eqref{CP16'} (respectively, \eqref{CP15'}) relates to the self-similar solution \eqref{CP9'} (respectively, to its interface, as in \eqref{CP10'}). If ~$\beta(p-1)=1$~we then have explicit values of $\zeta_*$~and ~$f_1(\rho)$ via \eqref{CP8'}, while in general we have lower and upper bounds via \eqref{CP11'}-\eqref{CP14'}. If~$u_0$~satisfies \eqref{CP3} with ~$\alpha=p/(p-1-\beta), C=C_*$,~then the small-time behavior of the interface and local solution depend on the terms smaller than ~$C_*(-x)^{p/(p-1-\beta)}$~in the expansion of ~$u_0~~~\text{as}~~x\rightarrow 0-$.\\

{\bf Region (3):} Let $b>0, 0<\beta<1, \alpha >p/(p-1-\beta)$~(here again we describe the results for the case ~$p=2$~as well). In this case the interface initially shrinks and
\begin{equation}\label{CP17'}\eta(t)\sim - \ell_*t^{1/\alpha(1-\beta)}~\text{as}~t\rightarrow 0+\end{equation}
where~$\ell_*=C^{-1/\alpha}(b(1-\beta))^{1/\alpha(1-\beta)}.$~For $\forall \ell>\ell_*$we have
\begin{equation}\label{CP18'}u(x,t)\sim\big [C^{1-\beta}(-x)_+^{\alpha(1-\beta)}-b(1-\beta)t\big]^{1/(1-\beta)}~\text{as}~t\rightarrow 0+\end{equation}
along the curve ~$x=\eta_l(t)=-lt^{1/\alpha(1-\beta)}.$ Hence, the interface initially coincides with that of the solution
\[\bar u(x,t)\sim\big [C^{1-\beta}(-x)_+^{\alpha(1-\beta)}-b(1-\beta)t\big]_+^{1/(1-\beta)}\]
to the problem 
\[\bar u_t+b\bar u^{\beta}=0,~~~~~~~~~~~\bar u(x,0)=C(-x)_+^{\alpha}.\]
Respective lower and upper estimations are given in section 4 (see (4.16) and (4.19) below). \\

{\bf Region (4):} In this case the interface initially has a waiting time. We divide the results into four different subcases (see Figure 1).\\
%\end{itemize}
%\begin{itemize}

{\bf(4a)} Let ~$\beta=1, \alpha =p/(p-2).$~This case reduces to the case ~$b=0$~by a simple transformation (see section 3). If~$u_0$~is defined by \eqref{CP4}, then the unique solution to \eqref{CP1}, \eqref{CP4} is
\begin{equation}\label{CP19'}u_C(x,t)=C(-x)_+^{p/(p-2)}\text{exp}(-bt)\big[1-(C/\bar C)^{p-2}b^{-1}(1-\text{exp}(-b(p-2)t))\big]^{1/(p-2)}\end{equation}
for~$x\in \mathbb{R},~~t\in [0,T),$where
\[T=+\infty\quad \text{if}~~b\geq (C/\bar C)^{p-2},\]
\[~~~~~~~~~~~~~~~~~~~~~~~~~~\quad \quad~~~~~~~~~~~~~T=(b(2-p))^{-1}\text{ln}[1-b(\bar C/C)^{p-2}],\quad \text{if}~~~-\infty<b<(C/\bar C)^{p-2}\]
\[~~~~~~~~~\bar C=\big[(p-2)^p/(2(p-1)p^{p-1})\big]^{1/(p-2)}\]
If~$u_0$~satisfies \eqref{CP3}, then lower and upper estimations are given by~$u_{C\pm \epsilon}.$\\

{\bf(4b)} Let~$\beta=1, \alpha>p/(p-2).$ Then for $\forall \epsilon>0 \  \exists x_{\epsilon}<0$ and $\delta_{\epsilon}>0$ such that
\begin{equation}\label{CP20'}(C-\epsilon)(-x)_+^{\alpha}\text{exp}(-bt)\leq u(x,t)\leq (C+\epsilon)(-x)_+^{\alpha}\text{exp}(-bt)\end{equation}
\begin{equation*}
\times\big[1-\epsilon b^{-1}(p-2)^{-p}\big(1-\text{exp}(-b(p-2)t)\big)\big]^{1/2-p},~~~~~~x>x_{\epsilon},~~0\leq t\leq \delta_{\epsilon},\end{equation*}\\

{\bf(4c)} Let $1<\beta<p-1,~\alpha\geq p/(p-1-\beta)$. Then for $\forall \epsilon>0 \ \exists x_{\epsilon}<0$ and $\delta_{\epsilon}>0$ such that
\begin{equation}\label{CP21'}g_{-\epsilon}(x,t) \leq u(x,t) \leq g_{\epsilon}(x,t), ~~x\geq x_{\epsilon},~~~0\leq t \leq \delta_{\epsilon}\end{equation}
where
\[g_{\epsilon}(x,t)=\begin{cases}
[(C+\epsilon)^{1-\beta}|x|^{\alpha(1-\beta)}+b(\beta-1)(1-d_{\epsilon})t]^{1/(1-\beta)},~~~~~x_{\epsilon}\leq x <0,\\
0, \qquad \qquad \qquad \qquad ~~~~~~~~~~~~~~~~~~~~~~~~~~~~~~~~~~~~~~~~~~~~~~~x\geq 0,
\end{cases}\]
\[d_{\epsilon}=\begin{cases}
\epsilon~\text{sign}~b~~~~~~~~~~~~~~~~~~~~~~~~~~~~~~~~~~\text{if} ~\alpha>p/(p-1-\beta),\\
\Big(\big((C+\epsilon)/C_*\big)^{p-1-\beta}+\epsilon\Big)~\text{sign}~b~~~~\text{if}~\alpha=p/(p-1-\beta),
\end{cases}\]
and the constant $C_*$ is defined in (I(2)).\\

{\bf(4d)} Let either~$1<\beta<p-1,~p/(p-2)\leq \alpha<p/(p-1-\beta),$ or ~$\beta \geq p-1,~~\alpha \geq p/(p-2).$\\
If $\alpha=p/(p-2)$ then for $\forall \epsilon>0 \ \exists x_{\epsilon}<0$ and $\delta_{\epsilon}>0$ such that
\begin{equation}\label{CP22'}(C-\epsilon)(-x)_+^{p/(p-2)}(1-\gamma_{-\epsilon}t)^{1/(2-p)}\leq u\leq (C+\epsilon)(-x)_+^{p/(p-2)}(1-\gamma_{\epsilon}t)^{1/(2-p)}\end{equation}
where\begin{equation*}\gamma_{\epsilon}=\big[2(p-1)p^{p-1}(C+\epsilon)^{p-2}/(p-2)^{1-p}\big]+\epsilon\end{equation*}
While if ~$\alpha>p/(p-2)$ then for $\forall \epsilon>0 \ \exists x_{\epsilon}<0$ and $\delta_{\epsilon}>0$ such that \begin{equation}\label{CP23'}
(C-\epsilon)(-x)_+^{\alpha}\leq u \leq (C+\epsilon)(-x)_+^{\alpha}(1-\epsilon t)^{1/2-p)},\quad ~~~~~~x\geq x_{\epsilon},~~ 0\leq t\leq \delta_{\epsilon}.\end{equation}
%\end{itemize}

{\bf(II)}~~$b=0$. We divide this case into three subcases. \\
%REMARK 2.1.\emph{Using the same techniques as in the case ~$b \neq 0$~we derive some global estimations}(\emph{see}\eqref{CP27'}).\emph{The new %element here is that we have constructed lower and upper solutions to the coresponding nonlinear ODE for the function}~$f(\xi)$~in \eqref{CP25'}.
%\begin{itemize}

{\bf(1)} Let ~$p>2,~0<\alpha <p/(p-2).$~In this case the interface expands. First, assume that ~$u_0$~ is defined by \eqref{CP4}. Then if ~$\alpha = (p-1)/(p-2)$ the explicit solution to the problem \eqref{CP1}, \eqref{CP4} is
\begin{equation}\label{CP24'}u(x,t)=C(\xi_*t-x)_+^{(p-1)/(p-2)},\quad \quad \xi_*=C^{p-2}\Big (\frac{p-1}{p-2}\Big )^{p-1}.\end{equation}
If~$0<\alpha <p/(p-2),$~then the solution to \eqref{CP1}, \eqref{CP4} has the self-similar form \eqref{CP25'}
\begin{equation}\label{CP26'}\eta(t)=\xi_* t^{\frac{1}{p-\alpha(p-2)}},\quad ~~~~~~~~~~~0\leq t<+\infty,\end{equation}
where~$\xi_*$~and~$f$~satisfy \eqref{CP2'}, \eqref{selfsimilarODE}-\eqref{CP6'}. Moreover, we have
\begin{equation}\label{CP27'}C_4t^{\frac{\alpha}{p-\alpha(p-2)}}(\xi_3-\xi)_+^{\frac{p-1}{p-2}}\leq u\leq C_5t^{\frac{\alpha}{p-\alpha(p-2)}}(\xi_4-\xi)_+^{\frac{p-1}{p-2}},\end{equation}
\[0\leq x<+\infty,~~~~~0<t<+\infty,\]
where $\xi_3$ (respectively, $\xi_4$) is defined by the right-hand side of \eqref{CP5'}, where we replace $\xi''_*$ with $C^{\frac{p-2}{p-\alpha(p-2)}}\xi_1$ (respectively, with  $C^{\frac{p-2}{p-\alpha(p-2)}}\xi_2$) and 
\[C_4 =C^{p/(p-\alpha(p-2))}A_0\xi_3^{(p-1)/(2-p)},~~~~~~~~~C_5 =C^{p/(p-\alpha(p-2))}A_0\xi_4^{-(p-1)/(p-2)}.\]
In the particular case $\alpha=(p-1)(p-2)^{-1}$, when an explicit solution is given by \eqref{CP24'}, we have $\xi_3=\xi_4=\xi_*$ and both lower and upper estimations in \eqref{CP27'} lead to the explicit solution \eqref{CP24'}. In general, when $\alpha \neq (p-1)(p-2)^{-1}$ the precise value $\xi_*$ relates to the similarity ODE for $f(\xi)$ from \eqref{selfsimilarODE}, namely, $\xi_*=\text{sup}\{\xi: f(\xi)>0\}.$ If $u_0$ satsfies \eqref{CP3} with $(0<\alpha<p/(p-2))$, then \eqref{CP1'} and \eqref{CP3'} are valid. Lower and upper bounds for $f(\rho)$ follow from \eqref{CP27'}.\\

{\bf(2)} Let $p>2, \alpha =p/(p-2).$ In this case the interface initially has a waiting time. If $u_0$ is defined by \eqref{CP4}, then the explicit solution to \eqref{CP1}, \eqref{CP4} is
\begin{equation}\label{CP28'}u_C(x,t)=C(-x)_+^{\alpha}\big[1-(C/\bar{C})^{p-2}(p-2)t\big]^{1/(2-p)}\quad x\in \mathbb{R},~~0\leq t<T\end{equation}
where
\[T=(\bar C/C)^{p-2}(p-2)^{-1}\]
and the constant ~$\bar C$~is defined in (I(4)).\\
If ~$u_0$~ satisfies \eqref{CP3} with ~$\alpha=p/(p-2)$,~then lower and upper estimations are given by ~$u_{C\pm \epsilon}$.\\

{\bf(3)} Let ~$p>2, \alpha >p/(p-2).$~In this case also the interface initially remains stationary and for $\forall \epsilon>0 \ \exists x_{\epsilon}<0$ and  $\delta_{\epsilon}>0$~such that 
\begin{equation}\label{CP29'}(C-\epsilon)(-x)_+^{\alpha}\leq u \leq (C+\epsilon)(-x)_+^{\alpha}(1-\epsilon t)^{1/2-p)},\quad ~~~~~~x_{\epsilon}\leq x,~~ 0\leq t\leq \delta_{\epsilon}\end{equation}
%\end{itemize}

\section{ Preliminary results.} The mathematical theory of nonlinear p-Laplacian type degenerate parabolic equations is well developed. 
Throughout this paper we shall follow the definition of weak solutions and of supersolutions (or subsolutions) of the equation \eqref{CP1} in the following sense:
\begin{definition}
A measurable function $u\geq 0$ is a local weak solution (respectively sub- or supersolution) of \eqref{CP1} in $\mathbb{R} \times (0,T]$ if 
\begin{itemize}
\item  $u \in C_{loc}(0,T; L^2_{loc}(\mathbb{R}) \cap L^p_{loc}(0,T; W_{loc}^{1,p}(\mathbb{R})\cap L^{1+\beta}_{loc}(\mathbb{R}))$
\item For $\forall$ subinterval $[t_0, t_1] \subset (0,T]$ and for $\forall \mu_i \in C^1[t_0, t_1], \ i=1,2$ such that $\mu_1(t) < \mu_2(t)$ for $t\in [t_0, t_1]$
\begin{equation}\label{supersub}
\int_{\mu_1(t)}^{\mu_2(t)}u \phi dx \Big |_{t_0}^{t_1}+\int_{t_0}^{t_1}\int_{\mu_1(t)}^{\mu_2(t)}(-u\phi_t+|u_x|^{p-2}u_x\phi_x+bu^\beta \phi )dx dt = 0 \  ( resp.  \leq  or \geq \  0)
\end{equation}
where $\phi \in C_{x,t}^{2,1}(\overline{D})$ is an arbitrary function that equals zero when $x=\mu_i(t), t_0\leq t \leq t_1, i=1,2$, and
\[ D=\{(x,t): \mu_1(t)<x<\mu_2(t), t_0< t <t_1\}  \]
\end{itemize}
\end{definition}
The questions of existence and uniqueness of initial boundary value problems for \eqref{CP1}, comparison theorems, and regularity of weak solutions are known due to \cite{dibe-sv, dibe1, dibe2, kalashnikov2, kalashnikov3, tsutsumi, est-vazquez} etc. The proof of the following typical comparison result is standard.
\begin{lemma}\label{lemma1}
Let $g$ be a nonnegative and continuous function in $\overline{Q}$, where
\[ Q=\{(x,t): \eta_0(t)<x<+\infty, 0< t <T\leq +\infty\},  \]
$f$ is in $C_{x,t}^{2,1}$ in $Q$ outside a finite number of curves $x=\eta_j(t)$, which divide $Q$ into a finite number of subdomains $Q^j$, where $\eta_j \in C[0,T]$; for arbitrary $\delta >0$ and finite $\Delta\in (\delta, T]$ the function $\eta_j$ is absolutely continuous in $[\delta, \Delta ]$. Let $g$ satisfy the inequality 
\[  Lg\equiv g_t-\Big(|g_x|^{p-2}g_x\Big)_x+bg^{\beta} \geq 0, \ (\leq 0)  \]
at the points of $Q$, where $g \in C_{x,t}^{2,1}$. Assume also that the function $|g_x|^{p-2}g_x$ is continuous in $Q$ and $g\in L^{\infty}(Q\cap (t\leq T_1))$ for any finite $T_1\in (0,T]$. Then $g$ is a supersolution (subsolution) of \eqref{CP1}. If, in addition we have
\[ g\Big |_{x=\eta_0(t)} \geq (\leq) \ u\Big |_{x=\eta_0(t)}, \ g\Big |_{t=0} \geq (\leq) \ u\Big |_{t=0} \]
then
\[ g \geq (\leq) \ u, \quad\text{in} \ \   \overline{Q}  \]
\end{lemma}

Suppose that ~$b \geq 0$~and that ~$u_0$~may have unbounded growth as ~$|x|\rightarrow +\infty$.~It is well known that in this case some restriction must be imposed on the growth rate for existence, uniqueness and comparison results in the CP \eqref{CP1}, \eqref{CP2}. Optimal growth condition for the equation (\eqref{CP1} with ~$b=0, p>2$ was derived in \cite{dibe2, dibe1}. If initial data may be majorized by power law function \eqref{CP4}, then there exists a unique solution (with $T=+\infty$) and a comparison principle is valid if ~$0<\alpha<p/(p-2)$. If ~$\alpha=p/(p-2)$, then existence, uniqueness, and comparison results are valid only locally in time. In particular, from \cite{dibe2,dibe1} it follows that the unique explicit solution to \eqref{CP1}, \eqref{CP4} with ~$b=0, \alpha=p/(p-2), T= (\bar C/C)^{p-2}(p-2)^{-1}$ is $u_C(x,t)$ from \eqref{CP28'}.

If the function ~$u(x,t)$~ is a solution to CP \eqref{CP1}, \eqref{CP4} with ~$b=0$, then the function
\[\bar u(x,t)=\text{exp}(-bt)u(x,(b(2-p))^{-1}\big(\text{exp}(b(2-p)t)-1\big))\]
is a solution to \eqref{CP1} with~$ b\neq 0, \beta=1$. Hence, from the above mentioned result it follows that the unique solution to CP \eqref{CP1}, \eqref{CP4} with ~$p>2, b\neq 0, \beta =1, \alpha=p/(p-2)$ is the function~$\bar u_C(x,t)$ from \eqref{CP19'}.

We are not interested in necessary and sufficient conditions on the growth rate at infinity for existence, uniqueness, and comparison results for the CP \eqref{CP1}, \eqref{CP2} with $b>0, p>2, \beta>0$; for our purposes it is enough to mention that if $u_0$ may be majorized by the function \eqref{CP4} with ~$\alpha$ satisfying~$ 0<\alpha< p/(p-2)$, then the CP \eqref{CP1}, \eqref{CP2} with $b>0, p>2, \beta>0, T=+\infty$ has a unique solution and for this class of initial data a comparison principle is valid. This easily follows from the fact that the solution of the CP \eqref{CP1}, \eqref{CP2} with $b=0$ is a supersolution of the CP with $b>0$, and hence it becomes a global locally bounded uniform upper bound for the increasing sequence of approximating bounded solutions of the CP with $b>0$.

In the next four lemmas we apply rescaling to establish some preliminary estimations of the solution to CP. 
\begin{lemma}\label{lemma2} 
If $b=0$ and $p>2, 0 <\alpha < p/(p-2),$ then the solution $u$ of the CP \eqref{CP1}, \eqref{CP4} has a self-similar form \eqref{CP25'}, where the self-similarity function $f$ satisfies \eqref{24ab}. If $u_0$ satisfies \eqref{CP3}, then the solution to CP \eqref{CP1}, \eqref{CP2} satisfies \eqref{CP1'}-\eqref{CP3'}
\end{lemma}
\begin{lemma}\label{lemma3}
Let u be a solution to the CP\eqref{CP1}, \eqref{CP2} and ~$u_0$ satisfy \eqref{CP3}. Let one of the following conditions be valid:

(a)~~~$b>0,~~~0<\beta<1<p,~~~0<\alpha<p/(p-1-\beta).$

(b)~~~$b\neq0,~~~\beta\geq 1,~~~p>2, ~~~0<\alpha<p/(p-2).$\\
Then ~$u$~satisfies \eqref{CP3'}.
\end{lemma}
\begin{lemma}\label{lemma4}
Let u be a solution to the CP\eqref{CP1}, \eqref{CP4} with~$b>0,~0<\beta<1, ~p>2,~\alpha=p/(p-1-\beta).$ Then the solution $u$ has the self-similar form \eqref{CP9'}. If ~$C>C_*$~then~$f_1(0)=A_1$,~where~$A_1$~is a positive number depending on ~$p, \beta, C$~and ~$b$. If~$u_0$~satisfies \eqref{CP3} with~$\alpha=p/(p-1-\beta), C>C_*, $~then ~$u$~satisfies
\begin{equation}\label{CP1''}u(0,t)\sim A_1 t^{1/(1-\beta)}~\text{as}~t\rightarrow 0+.\end{equation}
\end{lemma}
\begin{lemma}\label{lemma5}
Let~$u$~be a solution to the CP \eqref{CP1}-\eqref{CP3} with~$b>0, 0<\beta<1, \alpha>p/(p-1-\beta)$.~Then for arbitrary~$\ell >\ell_*$ (see \eqref{CP17'}) the asymptotic formula \eqref{CP18'} is valid with $x=\eta_{\ell}(t)=-\ell t^{1/\alpha(1-\beta)}.$
\end{lemma}
\emph{Proof of Lemma~\ref{lemma2}}. If we consider a function 
\begin{equation}\label{CP2''}u_k(x,t)=ku(k^{-1/\alpha}x,k^{(\alpha(p-2)-p)/\alpha}),~~~~~~k>0,\end{equation}
it may easily be checked that this satisfies \eqref{CP1}, \eqref{CP4}. From \cite{dibe1, dibe2} it follows that under the condition of the lemma there exists a unique global solution to \eqref{CP1}, \eqref{CP4}. Therefore, we have 
\begin{equation}\label{CP3''}u(x,t)=ku(k^{-1/\alpha}x,k^{(\alpha(p-2)-p)/\alpha}),~~~~~~k>0.\end{equation}
If we choose~$k=t^{\alpha/(p-\alpha(p-2))},$~then \eqref{CP3''} implies \eqref{CP25'} for $u$ with~$f(\xi)=u(\xi,1).$ In fact, $f$ is a unique nonnegative and differentiable weak solution of the  boundary value problem 
\begin{equation}\begin{cases}
\big(|f'(\xi)|^{p-2}f'(\xi)\big)'+\frac{1}{p-\alpha (p-2)}\xi f'(\xi)-\frac{\alpha}{p-\alpha (p-2)}f(\xi)=0,~-\infty <\xi <+\infty\\
f(-\infty)\sim C(-\xi)^\alpha, \  f(+\infty)= 0
\end{cases}\label{selfsimilarODE1}\end{equation}
 and there exists an~$\xi_*>0$~such that~$f$ satisfies \eqref{selfsimilarODE}: it is positive and smooth for~$\xi<\xi_*$~and~$f=0$~for~$\xi\geq \xi_*(\cite{Barenblatt1}).$~Thus, \eqref{CP26'} is valid. To find the dependence of~$f$~on~$C$~we can again use scaling. Namely, let~$w$~be a solution of the CP \eqref{CP1}, \eqref{CP4} with~$C=1$.~Then it may be easily checked that forarbitrary~$k>0$~
\[u(x,t)=kw(C^{1/\alpha}k^{-1/\alpha}x,C^{p/\alpha}k^{(\alpha(p-2)-p)/\alpha}t).\]
By choosing~$k=(C^{p/\alpha}t)^{\alpha/(p-\alpha(p-2))}$~we then have
\begin{equation}\label{CP4''}u(x,t)=C^{\frac{p}{p-\alpha(p-2)}}w(C^{\frac{p-2}{\alpha(p-2)-p}}\xi,1)t^{\alpha/(p-\alpha(p-2))}.\end{equation}
From \eqref{CP4''} and \eqref{CP25'}, \eqref{24ab} and \eqref{CP2'} follow.

Now assume that~$u_0$~satisfies \eqref{CP3}. Then for $\forall$ sufficiently small~$\epsilon>0 \ \exists \ x_{\epsilon}<0$~such that
\begin{equation}\label{CP5''}(C-\epsilon /2)(-x)_+^{\alpha}\leq u_0(x) \leq (C+\epsilon /2)(-x)_+^{\alpha},~~~~~~~~~~x\geq x_{\epsilon}.\end{equation}
Let~$u_{\epsilon}(x,t)$~(respectively,~$u_{-\epsilon}(x,t)$)~be a solution to the CP \eqref{CP1}, \eqref{CP2} with initial data~$(C+\epsilon)(-x)_+^{\alpha}$~(respectively,~$(C-\epsilon)(-x)_+^{\alpha}$). Since the solution to the CP \eqref{CP1}, \eqref{CP2} is continous there exists a number~$\delta=\delta(\epsilon)>0$~such that
\begin{equation}\label{CP6''}u_{\epsilon}(x_{\epsilon},t)\geq u(x_{\epsilon},t),~~~~u_{-\epsilon}(x_{\epsilon},t)\leq u(x_{\epsilon},t)~\text{for}~0\leq t \leq \delta.\end{equation}
From \eqref{CP5''}, \eqref{CP6''}, and a comparison principle, it follows that
\begin{equation}\label{CP7''}u_{-\epsilon}\leq u \leq u_{\epsilon}~~\text{for}~~x\geq x_{\epsilon},~~~~~~~0\leq t\leq \delta.  \end{equation}
Obviously
\begin{equation}\label{CP8''}u_{\pm \epsilon}(\xi_{\rho}(t),t)=f(\rho;C\pm \epsilon)t^{\alpha/(p-\alpha(p-2))},~~~t\geq 0.\end{equation}
(Furthermore, we denote the right-hand side of\eqref{CP4a'} by~$f(\rho, C)$.) Now taking~$x=\xi_{\rho}(t)$~in \eqref{CP7''}, after multiplying to~$t^{-\alpha/(p-\alpha(p-2))}$~and passing to the limit, first as~$t\rightarrow 0$~and then as~$\epsilon\rightarrow 0,$~
we can easily derive \eqref{CP3'}. Similarly, from \eqref{CP7''},~\eqref{CP26'}, and \eqref{CP2'}, \eqref{CP1'} easily follows. The lemma is proved.

\emph{Proof of Lemma~\ref{lemma3}}. As in the previous proof,  \eqref{CP5''}, \eqref{CP6''} and \eqref{CP7''} follow from \eqref{CP3}. Let the conditions of one of the cases (a) or (b) with~$b>0$~be valid. Then from the results mentioned earlier it follows that the existence, uniqueness, and comparison results of the CP \eqref{CP1}, \eqref{CP2} with~$u_0=(C\pm \epsilon)(-x)_+^{\alpha},~~T=+\infty$~hold. Now if we rescale
\begin{equation}\label{CP9''}u_k^{\pm \epsilon}(x,t)=ku_{\pm \epsilon}\big(k^{-1/\alpha}x, k^{(\alpha(p-2)-p)/\alpha}t\big),~~~~~~k>0  \end{equation}
then~$u_k^{\pm}(x,t)$~satisfies the following problem:
\begin{subequations}\label{3.12}
\begin{equation}\label{CP10a''}u_t-(|u_x|^{p-2}u_x)_x+bk^{(\alpha(p-1-\beta)-p)/\alpha}u^{\beta}=0,~~~~~~x\in \mathbb{R},~~~~t>0,  \end{equation}
\begin{equation}\label{CP10b''}u(x,0)=(C\pm \epsilon)(-x)_+^{\alpha},~~~~~~x\in \mathbb{R}.  \end{equation}
\end{subequations}
There exists a unique solution to CP \eqref{3.12}, which also obeys a comparison principle. Since~$\alpha(p-1-\beta)-p<0,$~by using a comparison principle in Lemma~\ref{lemma2} it follows that
\begin{equation}\label{monoton}
u_{k_1}^{\pm \epsilon}(x,t)\leq u_{k_2}^{\pm \epsilon}(x,t)\leq \cdots \leq v_{\pm}(x,t),~~~~~x\in \mathbb{R},~~~t\geq 0; \ \quad\text{if} \ k_1<k_2,
\end{equation}
where~$v_{\pm \epsilon}$~is a solution to CP  \eqref{CP1},  \eqref{CP2} with~$b=0,~u_0=(C\pm \epsilon)(-x)_+^{\alpha},\;T=+\infty$. From the results of \cite{dibe1, tsutsumi} it follows that the sequence of nonnegative and locally bounded solutions $\{u_{k}^{\pm \epsilon}\}$ is locally uniformly H\"{o}lder continuous, and weakly precompact in $W_{loc}^{1,p}(\mathbb{R}\times(0,T))$. Since $\alpha(p-1-\beta)-p<0,$ passing to limit as $k\to +\infty$, from \eqref{supersub} it follows that the limit function is a solution of the CP \eqref{CP1},  \eqref{CP2} with~$b=0,~u_0=(C\pm \epsilon)(-x)_+^{\alpha},\;T=+\infty$. Due to uniqueness we have
\begin{equation}\label{CP11''}\underset{k\rightarrow+\infty}\lim u_k^{\pm \epsilon}(x,t)=v_{\pm}(x,t),~~~~~x\in \mathbb{R},~~~t\geq 0,  \end{equation}
Hence,~$v_{\pm \epsilon}$~satisfies  \eqref{CP8''}. If we now take~$x=\xi_{\rho}(t),$~where~$\rho$~is an arbitrary fixed number satisfying~$\rho<\xi_*,$~then from  \eqref{CP11''} it follows that
\begin{equation}\label{CP12''}\underset{k\rightarrow+\infty}\lim ku_{\pm \epsilon}\big(k^{-1/\alpha}\xi_{\rho}(t), k^{(\alpha(p-2)-p)/\alpha}t\big)=f(\rho;C\pm \epsilon)t^{\alpha(p-\alpha(p-2))},~~~~~t>0.  \end{equation}
If we take~$\tau= k^{(\alpha(p-2)-p)/\alpha}t,$~then  \eqref{CP12''} imples\
\begin{equation}\label{CP13''}u_{\pm \epsilon}(\xi_{\rho}(\tau),\tau)\sim f(\rho;C\pm \epsilon)\tau^{\alpha(p-\alpha(p-2))},~\text{as}~\tau\rightarrow 0+.  \end{equation}
As before,  \eqref{CP3'} follows from  \eqref{CP7''},  \eqref{CP13''}.

Now consider the case (b) with~$b<0$.~Suppose that~$u_{\pm \epsilon}$~is a solution of the Dirichlet problem\begin{subequations}
\begin{equation}\label{CP14a''}u_t-(|u_x|^{p-2}u_x)_x+bu^{\beta}=0,~~~|x|<|x_{\epsilon}|,\;0<t<\delta,  \end{equation}
\begin{equation}\label{CP14b''}u(x,0)=(C\pm \epsilon)(-x)_+^{\alpha},~~~|x|\leq |x_{\epsilon}|,  \end{equation}
\begin{equation}\label{CP14c''}u(x_{\epsilon},t)=(C\pm \epsilon)(-x_{\epsilon})^{\alpha},~~~~u(-x_{\epsilon},t)=0,\;~0\leq t\leq \delta.  \end{equation}
\end{subequations}
The function~$u_k^{\pm \epsilon}, $~defined as in  \eqref{CP9''}, satisfies the Dirichlet problem
\begin{subequations}
\begin{equation}\label{CP15a''}
u_t-(|u_x|^{p-2}u_x)_x+bk^{(\alpha(p-1-\beta)-p)/\alpha}u^{\beta}=0~\text{in}~D^k_{\epsilon},
\end{equation}
\begin{equation}\label{CP15b''}
u(k^{1/\alpha}x_{\epsilon},t)=k(C\pm \epsilon)(-x_\epsilon)^{\alpha},~~~u(-k^{1/\alpha}x_{\epsilon},t)=0,~~~~~~0\leq t\leq k^{(p-\alpha(p-2))/\alpha}\delta
\end{equation}
\begin{equation}\label{CP15c''}
u(x,0)=(C\pm \epsilon)(-x)_+^{\alpha},~~~~~|x|\leq k^{1/\alpha}|x_{\epsilon}|,\end{equation}
\end{subequations}
where
\[D^k_{\epsilon}=\{(x,t):  |x|<k^{1/\alpha}|x_{\epsilon}|,~~~0<t\leq k^{(p-\alpha(p-2))/\alpha}\delta\}.\]
There exists a number~$\delta>0$~(which does not depend on $k$) such that both \eqref{CP14a''}-\eqref{CP14c''} and \eqref{CP15a''}-\eqref{CP15c''} have a unique solution (see discussion preceding Lemma~\ref{lemma2}).
In view of finite speed of propagation a~$\delta=\delta(\epsilon)>0$~may be chosen such that
\begin{equation}\label{CP16''}u(-x_{\epsilon},t)=0,~~~~0\leq t \leq \delta.\end{equation}
Applying the comparison theorem, from \eqref{CP5''}, \eqref{CP6''} and \eqref{CP16''},\eqref{CP7''} follows for~$|x|\leq |x_{\epsilon}|,\;0\leq t\leq\delta.$

To prove the convergence of the sequences~$\{u_k^{\pm \epsilon}\}$~as~$k\rightarrow +\infty$, we need to prove uniform boundedness. Consider a function
\[g(x,t)=(C+1) (1+x^2)^{\frac{\alpha}{2}}(1-\nu t)^{\frac{1}{2-p}},~x\in \mathbb{R},~0\leq t\leq t_0=\frac{\nu ^{-1}}{2}\]
where
\[\nu=h_*+1,~h_*=h_*(\alpha;p)=\max_{x\in\mathbb{R}}h(x),\]
\[h(x)=(p-2)\alpha^{p-1}(C+1)^{p-2}(1+x^2)^{\frac{(\alpha-2)(p-1)-2-\alpha}{2}}x^{2}|x|^{p-2}\Big(\frac{1+x^2}{x^2}+(p-2)\frac{1+x^2}{|x|^{2}}+(\alpha-2)(p-1)\Big)\]
Then we have
\[L_{k}g\equiv g_t-\big(|g_x|^{p-2}g_x\big)_x+bk^{\frac{\alpha(p-\beta -1)-p}{\alpha}}g^{\beta}=(C+1)(p-2)^{-1}(1+x^2)^{\frac{\alpha}{2}}(1-\nu t)^{\frac{p-1}{2-p}}S~~\text{in}~D^k_{\epsilon},\]
\[S=\nu-h(x)+b(p-2)(C+1)^{\beta -1}k^{\frac{\alpha(p-\beta -1)-p}{\alpha}}(1+x^2)^{\frac{\alpha(\beta-1)}{2}}(1-\nu t)^{\frac{\beta +1-p}{2-p}},\]
and hence
\begin{equation}\label{CP17''}S\geq 1+R~\text{in}~D_{0\epsilon}^{k}=D_{\epsilon}^{k}\cap \{0<t\leq t_{0}\},\end{equation}
where
\[R=O\Big(k^{p-2-p/\alpha}\Big)~~~\text{uniformly for}~~(x,t)\in D_{0\epsilon}^{k}~~\text{as}~~k\rightarrow +\infty.\]
Moreover, we have for~~$0<\epsilon \ll 1$
\begin{subequations}
\begin{equation}\label{CP18a''} g(x,0)~~~~~\geq u_k^{\pm \epsilon}(x,0)~\text{for}~|x|\leq k^{1/\alpha}|x_{\epsilon}|,\end{equation}
\begin{equation}\label{CP18b''}g(\pm k^{1/\alpha}x_{\epsilon},t)\geq u_k^{\pm \epsilon}(\pm k^{1/\alpha}x_\epsilon,t)~\text{for}~0\leq t\leq t_0.\end{equation}
\end{subequations}
Hence, $\exists \ k_0=k_0(\alpha;p)$~such that for $\forall k\geq k_0$~the comparison theorem implies
\begin{equation}\label{CP19''}0 \leq u_k^{\pm \epsilon}(x,t)\leq g(x,t)~\text{in}~\bar D_{0\epsilon}^k.\end{equation}
Let~$ G$~ be an arbitrary fixed compact subset of 
\[P=\big\{(x,t):x\in \mathbb{R},~~~0<t\leq t_0\big\}.\]
We take~$k_0$~so large that~$G\subset D_{0\epsilon}^k$~for~$k\geq k_0.$~From \eqref{CP19''} it follows that the sequences~$\{u_k^{\pm \epsilon}\}$,~$k\geq k_0 $,~are uniformly bounded in~$G$.~As before, from the results of \cite{dibe1, tsutsumi} it follows that the sequence of nonnegative and locally bounded solutions $\{u_{k}^{\pm \epsilon}\}$ is locally uniformly H\"{o}lder continuous, and weakly precompact in $W_{loc}^{1,p}(\mathbb{R}\times(0,T))$. It follows that for some subsequence~$k'$~
\begin{equation}\label{CP20''}\underset{k'\rightarrow+\infty}\lim u_{k'}^{\pm \epsilon}(x,t)=v_{\pm \epsilon}(x,t),~~~~~~~(x,t)\in P.\end{equation}
Since $\alpha(p-1-\beta)-p<0,$ passing to limit as $k'\to +\infty$, from \eqref{supersub} for $u_{k'}^{\pm \epsilon}$ it follows that~$v_{\pm \epsilon}$~is a solution to the CP \eqref{CP1}, \eqref{CP2} with~$b=0, T=t_0, u_0=(C\pm \epsilon)(-x)_+^{\alpha}.$~As before, from \eqref{CP8''}, \eqref{CP12''}, \eqref{CP13''} and \eqref{CP7''}, the required estimation \eqref{CP3'} follows. The lemma is proved.

The first assertion of Lemma~\ref{lemma4}  has been proved in \cite{Lidumu} for the case~$p>2.$ If~$u_0$~satisfies \eqref{CP3}, the estimation \eqref{CP1''} may be proved exactly as estimation \eqref{CP3'} was proved in Lemma~\ref{lemma2}. 

\emph{Proof of Lemma~\ref{lemma5}}. Asymptotic behaviour  \eqref{CP3} imply \eqref{CP5''} and \eqref{CP6''}. Assume that that~$v_{\pm \epsilon}$ solves  the problem 
\[v_t-(|v_x|^{p-2}v_x)_x+bv^{\beta}=0,\quad \quad |x|<|x_{\epsilon}|,~~0<t\leq \delta,\]\begin{equation*}
v(x,0)=(C\pm \epsilon)(-x)_+^{\alpha},\quad \quad ~~\quad|x|\leq |x_{\epsilon}|,\end{equation*}
\begin{equation*}
v(x_{\epsilon},t)=(C\pm \epsilon)(-x_{\epsilon})_+^{\alpha},~~ v(-x_{\epsilon},t)=u(-x_{\epsilon},t), ~\quad 0\leq t\leq \delta.\end{equation*}
According to comparison result from \eqref{CP5''} and \eqref{CP6''}, \eqref{CP7''} follows for~$|x|\leq |x_{\epsilon}|,~~0\leq t \leq \delta.$ If we rescale
\[u_k^{\pm \epsilon}(x,t)=ku_{\pm \epsilon}(k^{-\frac{1}{\alpha}}x,k^{\beta-1}t),~~~~~k>0,\]
then~$u_k^{\pm \epsilon}$~satisfies the Dirichlet problem
\begin{equation*}
v_t-k^{\frac{p-\alpha(p-1-\beta)}{\alpha}}\big(|v_x|^{p-2}v_x\big)_x+bv^{\beta}=0 ~\text{in}~ E^k_{\epsilon}\end{equation*}
\begin{equation*}
v(x,0)=(C\pm \epsilon)(-x)_+^{\alpha},~~|x|\leq k^{\frac{1}{\alpha}}|x_{\epsilon}|\end{equation*}
\begin{equation*}
v(k^{\frac{1}{\alpha}}x_{\epsilon},t)=k(C\pm \epsilon)(-x_{\epsilon})_+^{\alpha},\quad v(-k^{\frac{1}{\alpha}}x_{\epsilon},t)=ku(-x_{\epsilon},k^{\beta-1}t),~0\leq t\leq k^{1-\beta}\delta\end{equation*}
where
\[E^k_{\epsilon}=\big\{|x|<k^{\frac{1}{\alpha}}|x_{\epsilon}|,~0<t\leq k^{1-\beta}\delta\big\}.\]
The goal is to in prove the convergence of the sequence~$\{u_k^{\pm \epsilon}\}$~as~$k\rightarrow +\infty$. To establish uniform bound consider $g(x,t)=(C+1)(1+x^2)^{\alpha/2}$~exp~$t$. We have
\begin{equation}\label{CP21''} \tilde{\text{L}}_kg \equiv g_t-k^{\frac{p-\alpha(p-1-\beta)}{\alpha}}\big(|g_x|^{p-2}g_x\big)_x+bg^{\beta}\geq g\Big[1-k^{\frac{p-\alpha(p--1-\beta)}{\alpha}}\alpha^{p-1}(C+1)^{p-2}e^{t(p-2)}\end{equation}
\begin{equation*}
\times(1+x^2)^{\frac{(\alpha-2)(p-1)-2-\alpha}{2}}x^2|x|^{p-2}\Big(\frac{1+x^2}{x^2}+(p-2)\frac{1+x^2}{|x|^{2}}+(\alpha-2)(p-1)\Big)\Big]~\text{in}~E_{\epsilon}^k.\end{equation*}
Let~$t_0>0$~be fixed and let~$E_{0\epsilon}^{k}=E_{\epsilon}^k\cap \{(x,t):~~0<t\leq t_0\}.$ ~From \eqref{CP21''} it follows that 
\[\tilde{L}_kg \geq (1+R)~~~~\text{in} ~~E_0^k,\]
where
\[R=O(k^{\theta})~\text{uniformly for}~(x,t)\in E_{0\epsilon}^k~~\text{as}~~k\rightarrow +\infty\]
\[\theta=\big(p-\alpha(p-1-\beta)/\alpha\big)~~~\text{if}~~\alpha<p/(p-2),\]
\begin{equation*}
\theta=\beta-1, \quad \quad \quad ~~\text{if}~ \alpha\geq p/(p-2).\end{equation*}
We have for~$0<\epsilon \ll{1}$~that
\[ g(x,0)=u_k^{\pm \epsilon}(x,0),~~~~\text{for}~~~|x|\leq k^{1/\alpha}|x_{\epsilon}|,\]
and
\[u_k^{\pm \epsilon}(-k^{\frac{1}{\alpha}}x_{\epsilon},t)=o(k),~0\leq t\leq t_0~\text{as}~ k\rightarrow \infty,\]
\[g(\pm k^{\frac{1}{\alpha}}x_{\epsilon},t)\geq u_k^{\pm \epsilon}(\pm k^{\frac{1}{\alpha}}x_{\epsilon},t),~~\text{for}~~~0\leq t\leq t_0\]
if~$k$~is chosen large enough. Therefore, the comparison principle implies \eqref{CP19''} in~$\bar E^k_{0\epsilon},$~where the respective functions~$u_k^{\pm \epsilon}$~and~$g$~apply in the context of the this proof. ~As before, from the interior regularity results (\cite{dibe1, tsutsumi}) it follows that the sequence of nonnegative and locally bounded solutions $\{u_{k}^{\pm \epsilon}\}$ is locally uniformly H\"{o}lder continuous, and weakly precompact in $W_{loc}^{1,p}(\mathbb{R}\times(0,T))$. It follows that for some subsequence~$k'$~, \eqref{CP20''} is valid. Since $ \alpha>p/(p-1-\beta)$, it follows that the limit functions $v_{\pm \epsilon}$ are solutions to the problem

\[v_t+bv^{\beta}=0,~x\in \mathbb{R},~0<t\leq t_0;\quad v(x,0)=(C\pm \epsilon)(-x)_+^{\alpha},~x\in \mathbb{R},\]
i.e.,
\[v_{\pm \epsilon}(x,t)=\Big[(C\pm \epsilon)^{1-\beta}(-x)_+^{\alpha(1-\beta)}-b(1-\beta)t\Big]^{\frac{1}{1-\beta}}_+.\]
Let~$l>l_*$~be an arbitrary number and~$\epsilon>0$~be chosen such that
\[(C-\epsilon)^{1-\beta}\ell^{\alpha(1-\beta)}>b(1-\beta).\]
If we now take ~$x=\eta_{\ell}(t)$~and $\tau =k^{\beta-1}t$,~it follows from \eqref{CP20''} that
\begin{equation}\label{CP22''}u_{\pm \epsilon}(\eta_{\ell}(\tau),\tau)\sim \Big[(C\pm \epsilon)^{1-\beta}\ell^{\alpha(1-\beta)}-b(1-\beta)\Big]^{\frac{1}{1-\beta}}\tau^{\frac{1}{1-\beta}}f~\text{as}~\tau\rightarrow 0+.\end{equation}
Since $\epsilon>0$ is arbitrary, From \eqref{CP7''} and \eqref{CP22''}, \eqref{CP18'} follows. The lemma is proved.

%{\bf Remark.}  \emph{Lemma~\ref{lemma5} is true also if $\beta<p-1\leq 1,$ the proof completely coinciding with the one given. Note that in this case %$\theta=(p-\alpha(p-1-\beta))/\alpha$ if~$\beta<p-1\leq 1.$} 

\section{Proofs of the main results.} In this section we prove the main results described in section 2.

(I)~$b\neq 0$~and~$p>2$.

(1) Assume~$\alpha<p/(p-1-\text{min}\{1,\beta\})$~The formula  \eqref{CP3'} follows from Lemma~\ref{lemma2}. Since~$\rho$ is arbitrary, it implies
\begin{equation}\label{CP1'''}\underset{t\rightarrow 0+}\lim \text{inf}~\eta(t)t^{1/(\alpha(p-2)-p)}\geq \xi_*.\end{equation}
Take an arbitrary sufficiently small number~$\epsilon>0$.~Let~$u_{\epsilon}$~be a solution of the CP  \eqref{CP1},  \eqref{CP4} with~$b=0$~and with $C$ replaced by $C+\epsilon$. As before, the second inequality of  \eqref{CP5''} and the first inequality of  \eqref{CP6''} follow from  \eqref{CP3}. Suppsoe that $b>0$. In this case, $u_{\epsilon}$ is a supersolution of  \eqref{CP1}. From  \eqref{CP5''},  \eqref{CP6''}, and a comparison principle, the second inequality of  \eqref{CP7''} follows. By Lemma 3.1 we then have
\[\eta(t)\leq(C+\epsilon)^\frac{2-p}{\alpha(p-2)-p}\xi'_*t^{1/(p-\alpha(p-2))},~~~0\leq t\leq\delta,\]
and hence 
\begin{equation}\label{CP2'''}\underset{t\rightarrow0^+}\lim\text{sup}~\eta(t)t^{\frac{1}{\alpha(p-2)-p}}\leq\xi_*.\end{equation}
Asssume now that $b<0$ and $\beta \geq 1$. The function
\[\bar u_{\epsilon}(x,t)=\text{exp}(-bt)u_{\epsilon}\Big(x,\frac{1}{b(2-p)}\big[\text{exp}(b(2-p)t)-1\big]\Big)\]
is a solution to the CP  \eqref{CP1},  \eqref{CP4} with $\beta=1$ and with $C$ replaced by $C+\epsilon$. As before, from  \eqref{CP3} the first inequality of  \eqref{CP6''} follows, where we replace $u_{\epsilon}$ with $\bar u_{\epsilon}$. Choose $|x_{\epsilon}|$ and $\delta $ so small that 
\[\bar u_{\epsilon} <1~\text {in}~B=\big\{(x,t):x\geq x_{\epsilon},~0<t\leq \delta\big\}.\]
Obviously, $\bar u_{\epsilon}$ is a supersolution of  \eqref{CP1} in $B$. From  \eqref{CP5''},  \eqref{CP6''}, and a comparison principle, the second inequality of  \eqref{CP7''}, with $u_{\epsilon}$ replaced by $\bar u_{\epsilon},$ follows. Thus we have
\[\eta(t)\leq (C+\epsilon)^{\frac{2-p}{\alpha(p-2)-p}} \xi'_*\Big\{\big(b(2-p)\big)^{-1}\big[\text{exp}(b(2-p)t)-1\big]\Big\}^{1/(p-\alpha(p-2))},~~~0\leq t\leq \delta,\]
which again implies \eqref{CP2'''}. From \eqref{CP1'''} and \eqref{CP2'''}, \eqref{CP1'} follows. Finally, \eqref{CP5'}, \eqref{CP6'}, \eqref{CP7'} follow from \eqref{CP27'}, which will be proved later in this section.

(2)~$b>0, 0<\beta<1, p>2,\alpha=p/(p-1-\beta).$~\\
First, consider the global case of \eqref{CP4}. The problem \eqref{CP1}, \eqref{CP4} has a unique global solution and for this class of initial data a comparison principle is valid (\cite{dibe1,dibe2}).\\
If~$\beta(p-1)=1$~it may be easily checked that the explicit solution to \eqref{CP1}, \eqref{CP4} is given by \eqref{CP8'}.\\
Let~$\beta(p-1)\neq 1.$~The self-similar form \eqref{CP9'} follows from Lemma~\ref{lemma4}. Let~$C>C_*.$
Consider a function
\begin{equation}\label{CP3'''}g(x,t) =t^{1/(1-\beta)}f_1(\zeta),~~~~~\zeta=xt^{-\frac{p-1-\beta}{p(1-\beta)}}.\end{equation}
we then have 
\begin{subequations}
\begin{equation}\label{CP4a'''}\text{L}g=t^{\frac{\beta}{1-\beta}}\mathcal{L}^0f_1,
\end{equation}
\begin{equation}\label{CP4b'''}\mathcal{L}^0f_1=\frac{1}{1-\beta}f_1-\big(|f'_1|^{p-2}f'_1\big)'-\frac{p-1-\beta}{p(1-\beta)}\zeta f'_1+bf_1^{\beta}.
\end{equation}
\end{subequations}
Choose as a function~$f_1$
\[f_1(\zeta)=C_0(\zeta_0-\zeta)_+^{\gamma_0},~~~~0<\zeta<+\infty,\]
where $~C_0,\zeta_0,\gamma_0 ~$are some positive constants. Taking~$\gamma_0=p/(p-1-\beta)$, from \eqref{CP4b'''} we have
\begin{equation}\label{CP5'''}\mathcal{L}^0f_1=bC_0^{\beta}(\zeta_0-\zeta)_+^{\frac{p\beta}{p-1-\beta}}\Big\{1-\Big(\frac{C_0}{C_*}\Big)^{p-1-\beta}+ \frac{C_0^{1-\beta}}{b(1-\beta)}\zeta_0(\zeta_0-\zeta)_+^{\frac{\beta(1-p)+1}{p-1-\beta}}\Big\}\end{equation}
To prove an upper estimation we take $C_0=C_2, \zeta_0=\zeta_2$~(see Appendix). If~$\beta(p-1)>1$,~then we have 
\[\mathcal{L}^0f_1\geq bC_2^{\beta}(\zeta_2-\zeta)_+^{\frac{p\beta}{p-1-\beta}}\Big\{1-\Big(\frac{C_2}{C_*}\big)^{p-1-\beta}+ \frac{C_2^{1-\beta}}{b(1-\beta)}\zeta_2^{\frac{p(1-\beta)}{p-1-\beta}}\Big\}=0,~ \text{for} ~0\leq \zeta \leq \zeta_2,\]
while if~$\beta(p-1)<1$,~we have 
\[\mathcal{L}^0f_1\geq bC_2^{\beta}(\zeta_2-\zeta)_+^{\frac{p\beta}{p-1-\beta}}\Big\{1-\Big(\frac{C_2}{C_*}\Big)^{p-1-\beta}\Big\}=0,~ \text{for} ~0\leq \zeta \leq \zeta_2.\]
From \eqref{CP4a'''}it follows that
\begin{subequations}
\begin{equation}\label{CP6a'''}\text{L}g\geq 0\quad \text{for}~0<x<\zeta_2 t^{\frac{p-1-\beta}{p(1-\beta)}},\;0<t<+\infty, 
\end{equation}
\begin{equation}\label{CP6b'''}\text{L}g= 0\quad \text{for}~x>\zeta_2 t^{\frac{p-1-\beta}{p(1-\beta)}},\;0<t<+\infty.\end{equation}
\end{subequations}
Lemma~\ref{lemma1} implies that~$g$~that $g$ is a supersolution of \eqref{CP1} in$ \{(x,t): x>0, t>0\}$. Since 
\begin{subequations}
\begin{equation}\label{CP7a'''}g(x,0)=u(x,0)=0\quad \text{for}~0\leq x<+\infty.
\end{equation}
\begin{equation}\label{CP7b'''}g(0,t)=u(0,t)\quad \text{for}~0\leq x<+\infty.
 \end{equation}
\end{subequations}
the right-hand side of \eqref{CP11'} follows.
If $\beta(p-1)<1$ then to prove the lower estimation we take $C_0=C_1, \zeta_0=\zeta_1,~\gamma_0=p/(p-1-\beta).$ Then from \eqref{CP5'''} we derive
\[\mathcal{L}^0f_1\leq bC_1^{\beta}(\zeta_1-\zeta)^{\frac{p\beta}{p-1-\beta}}\Big\{1-\big(\frac{C_1}{C_*}\big)^{p-1-\beta}+\frac{C_1^{1-\beta}}{b(1-\beta)}\zeta_1^{\frac{p(1-\beta)}{p-1-\beta}}\Big\}=0~\text{for}~0\leq\zeta \leq \zeta_1,\]
and from \eqref{CP4a'''} it follows that
\begin{subequations}
\begin{equation}\label{CP8a'''}\text{L}g\leq 0\quad \text{for}~0<x<\zeta_1 t^{\frac{p-1-\beta}{p(1-\beta)}},\;0<t<+\infty, \end{equation}
\begin{equation}\label{CP8b'''}\text{L}g= 0\quad \text{for} ~~x>\zeta_1 t^{\frac{p-1-\beta}{p(1-\beta)}},\;0<t<+\infty. \end{equation}
\end{subequations}
As before, from Lemma~\ref{lemma1} and \eqref{CP7a'''},\eqref{CP7b'''} the left-hand side of  \eqref{CP11'} follows.\\
If $\beta(p-1)>1,$ then to prove the lower estimation we take $C_0=C_1, \zeta_0=\zeta_1,~\gamma_0=(p-1)/(p-2).$ Then from  \eqref{CP4b'''} we have 
\[\mathcal{L}^0f_1=C_1(1-\beta)^{-1}(\zeta_1-\zeta)^{\frac{1}{p-2}}\Big\{\zeta_1-\Big(\frac{\beta(p-1)-1}{p(p-2)}\Big) \zeta-(1-\beta)C_1^{p-2}\big(\frac{p-1}{p-2}\big)^{p}+b(1-\beta)C_1^{\beta -1}(\zeta_1-\zeta)^{ \frac{\beta (p-1)-1}{p-2}}\Big\}\]
\[\leq C_1(1-\beta)^{-1}(\zeta_1-\zeta)^{\frac{1}{p-2}}\Big\{\zeta_1-C_1^{p-2}\frac{(1-\beta)(p-1)^p}{(p-2)^p}+b(1-\beta)C_1^{\beta -1}\zeta_1^{ \frac{\beta (p-1)-1}{p-2}}\Big\}=0~~~\text{for}~0<\zeta<\zeta_1,\]
which again implies \eqref{CP8a'''},\eqref{CP8b'''}. From Lemma~\ref{lemma1}, the left-hand side of  \eqref{CP11'} follows.

By applying the same analysis it may easily be checked that the alternative upper estimation is valied if $C_0=\bar C_2, \zeta_0=\bar \zeta_2, \gamma_0=(p-1)/(p-2).$

Let $\beta(p-1)>1$ and $0<C<C_*.$ Consider a function
\[g(x,t)=\big[C^{1-\beta}(-x)_+^{\frac{p(1-\beta)}{p-1-\beta}}-b(1-\beta)(1-\gamma)t\big]_+^{\frac{1}{1-\beta}},~~~~x\in\mathbb{R},~t>0,\]
where $\gamma\in [0,1).$ Let us estimate L$g$ in 
\[M=\{(x,t):-\infty<x<\mu_{\gamma}(t),~t>0\}, \ \ \ 
\mu_{\gamma}(t)=-[b(1-\beta)(1-\gamma)C^{\beta-1}t]^{\frac{p-1-\beta}{p(1-\beta)}}.\]
we have
\begin{subequations}
\begin{equation*}
Lg=bg^{\beta}S, \ S=\gamma-p^{p-1}(\beta(1-p)+1)(p-1)b^{-1}(p-1-\beta)^{-p}C^{p-1-\beta}\Big[1-\frac{b(1-\beta)(1-\gamma)t}{C^{1-\beta}(-x)_+^{\frac{p(1-\beta)}{p-1-\beta}}}\Big)\Big]^{\frac{\beta(p-2)}{1-\beta}} \end{equation*}
\begin{equation}\label{CP9a'''}-p^p\beta(p-1)b^{-1}(p-1-\beta)^{-p}C^{p-1-\beta}\Big[1-\frac{b(1-\beta)(1-\gamma)t}{C^{1-\beta}(-x)_+^{\frac{p(1-\beta)}{p-1-\beta}}}\Big]^{\frac{\beta(p-1)-1}{1-\beta}}. \end{equation}
Hence
\begin{equation}\label{CP9b'''}S|_{t=0}=\gamma-\Big(\frac{C}{C_*}\Big)^{p-1-\beta}, \quad S|_{x=\mu_{\gamma}(t)}=\gamma. 
\end{equation}
\end{subequations}
Moreover
\[S_t=\frac{p^{p-1}(p-1)(1-\gamma)C^{p-2}}{(p-1-\beta)^p}(-x)_+^{\frac{p(\beta-1)}{p-1-\beta}}\Big[1-C^{\beta-1}(-x)_+^{\frac{p(\beta-1)}{p-1-\beta}}b(1-\beta)(1-\gamma)t\Big]^{\frac{p\beta-2}{1-\beta}}\times\]
\[\times\Big[(\beta(p-1)-1)\beta(p-2)C^{\beta-1}b(1-\beta)(-x)_+^{-\frac{p(1-\beta)}{p-1-\beta}}(1-\gamma)t+(\beta(p-1)-1)(2\beta)\Big]\geq 0\quad \text{in}~M.\]
Thus
\[\gamma-\Big(\frac{C}{C_*}\Big)^{p-1-\beta}\leq S\leq\gamma \quad \text{in} ~M,\]
If we take ~$\gamma=\Big(\frac{C}{C_*}\Big)^{p-1-\beta}$ (respectively, $\gamma=0$), then we have\begin{subequations}
\begin{equation}\label{CP10a'''}
\text{L}g\geq 0 (\text{respectively}, \text{L}g\leq 0)~\text{in} ~M,\end{equation}
\begin{equation}\label{CP10b'''}\text{L}g= 0~\text{for}~x>\mu_{\gamma}(t),\;t>0\end{equation} 
\end{subequations}
and the estimation \eqref{CP13'} follows from the Lemma~\ref{lemma1}.

Let $\beta(p-1)<1$ and $0<C<C_*.$ First, we can establish the following rough estimation:
\begin{equation}\label{CP11'''}
\Big[C^{1-\beta}(-x)_+^{\frac{p(1-\beta)}{p-1-\beta}}-b(1-\beta)\big(1-\Big(\frac{C}{C_*}\Big)^{p-1-\beta}\big)t\Big]_+^{\frac{1}{1-\beta}}\leq u(x,t)
\leq C(-x)_+^{\frac{p}{p-1-\beta}} \ x\in \mathbb{R}, 0\leq t<+\infty.\end{equation}
To prove the left-hand side we consider the function $g$ as in the case when ~$\beta(p-1)>1$~with $\gamma = \big(C/C_*\big)^{p-1-\beta}.$~As before, we then derive \eqref{CP9a'''} and, since 
\[S_t\leq 0~~~~\text{in}~~~M,\]
we have ~$S\leq 0~\text{in}~ M.$ Hence, \eqref{CP10a'''},\eqref{CP10b'''} are valid with reversed inequality. As before, from Lemma~\ref{lemma1} the left-hand side of \eqref{CP11'''} follows. Since
\[\text{L}u_0=bu_0^{\beta}\big(1-\Big(C/C_*\Big)^{p-1-\beta}\big)\geq 0~~ ~\text{for}~~~x\in \mathbb{R},~t\geq0,\]
the second inequality in \eqref{CP11'''} follows.
Using \eqref{CP11'''}, we can now establish a more accurate estimation \eqref{CP14'}. Consider a function
\[g(x,t)=C_0(-\zeta_0 t^{\frac{p-1-\beta}{p(1-\beta)}}-x)^{\frac{p}{p-1-\beta}}_+~~~\text{in}~~~ G_{\ell},\]
\[G_{\ell}=\{(x,t):\zeta(t)=-\ell t^{\frac{p-1-\beta}{p(1-\beta)}}<x<+\infty,~0<t<+\infty\},\]
where,~$C_0>0,~\zeta_0>0,~\ell>\zeta_0$~are some constants. Calculating L$g$ in 
\[G^+_\ell=\{(x,t);\zeta(t)<x<-\zeta_0 t^{\frac{p-1-\beta}{p(1-\beta)}},~0<t<+\infty\},\]
we have
\[\text{L}g=bg^{\beta}S,\quad S=1-\Big(C_0/C_*\Big)^{p-1-\beta}-(b(1-\beta))^{-1}C_0^{1-\beta}\zeta_0 t^{\frac{\beta(p-1)-1}{p(1-\beta)}}\]
\begin{equation}\label{CP12'''}\times(-\zeta_0 t^{\frac{p-1-\beta}{p(1-\beta)}}-x)^{\frac{\beta(1-p)+1}{p-1-\beta}}.\end{equation}
Hence, if we take $C_0=C_*,$ then
\begin{equation}\label{CP13'''}\text{L}g\leq 0~\text{ in}~ G_{\ell}^+;~\text{L}g=0~\text{in}~G_{\ell}\backslash \bar G_{\ell}^+.\end{equation}
To obtain a lower estimation we now choose $\zeta_0=\zeta_3,~ \ell=\ell_0$ (see Appendix). Using \eqref{CP11'''}, we have
\[g(\zeta(t),t)=C_*(\ell_0-\zeta_3)^{\frac{p}{p-1-\beta}}t^{\frac{1}{1-\beta}}=\Big(b(1-\beta)\theta_*t\Big)^{\frac{1}{1-\beta}}\]
\begin{subequations}
\begin{equation}\label{CP14a'''}=\Big[C^{1-\beta}\ell_0^{\frac{p(1-\beta)}{p-1-\beta}}-b(1-\beta)\big(1-\Big(C/C_*\Big)^{p-1-\beta}\big)\Big]^{\frac{1}{1-\beta}}t^{\frac{1}{1-\beta}}
\leq u(\zeta(t),t),~t\geq0,\end{equation}
\begin{equation}\label{CP14b'''}g(x,0)=u(x,0)=0, ~0\leq x\leq x_0\end{equation}
\begin{equation}\label{CP14c'''}g(x_0,t)=u(x_0,t)=0, ~t\geq 0\end{equation}
\end{subequations}
where $x_0>0$ is an arbitrary fixed number. By using \eqref{CP13'''}, \eqref{CP14a'''}-\eqref{CP14c'''}, we can apply Lemma~\ref{lemma1} in
\[G'_{\ell_0}=G_{\ell_0}\cap \Big\{x<x_0\Big\}.\]
Since $x_0>0$ is arbitrary number the desired lower estimation from \eqref{CP14'} follows.

Let us now prove the right-hand side of \eqref{CP14'}. Since
\[S_x\geq 0, ~~~\text{for}~~~ \zeta(t)<x<-\zeta_0 t^{\frac{p-1-\beta}{p(1-\beta)}},~t>0,\]
from \eqref{CP12'''} it follows that
\[S\geq S|_{x=\zeta((t)}=1-\Big(C_0/C_*\Big)^{p-1-\beta}-\big(b(1-\beta)\big)^{-1}C_0^{1-\beta}\zeta_0 (\ell-\zeta_0 )^{\frac{\beta(1-p)+1}{p-1-\beta}}.\]
Taking now~$ C_0=C_3,~\zeta_0=\zeta_4,~\ell=\ell_1$ (see Appendix), we have
\[S|_{x=\zeta(t)}=0;\]
hence (by using \eqref{CP11'''})
\[\text{L}g\geq 0~~~\text{in}~~G^+_{\ell_1},~\text{L}g=0~~~\text{in}~~G_{\ell_1}\backslash \bar G^+_{\ell_1}\]
\[u(\zeta(t),t)\leq C\ell_1^{\frac{p}{p-1-\beta}}t^{\frac{1}{1-\beta}}=C_3(\ell_1-\zeta_4)^{\frac{p}{p-1-\beta}}t^{\frac{1}{1-\beta}}=g(\zeta(t),t),~t\geq 0,\]
and, for arbitrary $x_0>0,$ \eqref{CP14b'''} and \eqref{CP14c'''} are valid. As before, applying the Lemma~\ref{lemma1} in $G'_{\ell_1},$ we then derive the right-hand side of \eqref{CP14'}, since $x_0>0$ is arbitrary. From \eqref{CP11'}, \eqref{CP13'}, and \eqref{CP14'} it follows that 
\[\zeta_1t^{\frac{p-1-\beta}{p(1-\beta)}}\leq \eta(t)\leq \zeta_2t^{\frac{p-1-\beta}{p(1-\beta)}},~~~~0\leq t<+\infty,\]
where the constants $\zeta_1$ and $\zeta_2$ are chosen according to relevant estimations for $u$. 
If~$u_0$~satisfies \eqref{CP3} with~$\alpha=p/(p-1-\beta)$~and with~$C\neq C_*$,~then the asymptotic formulae \eqref{CP15'} and \eqref{CP16'}may be proved as the similar estimations \eqref{CP1'} and \eqref{CP3'} were in Lemma~\ref{lemma2}.

(3)~~Suppose that $b>0,~ 0<\beta<1,~\alpha>p/(p-1-\beta),~p>2$. 

Take an arbitrary sufficiently small number $\epsilon>0.$ From \eqref{CP3}, \eqref{CP5''} follows. Then consider a function
\begin{equation}\label{CP15'''}g_{\epsilon}(x,t)=\big[(C+\epsilon)^{1-\beta}(-x)_+^{\alpha(1-\beta)}-b(1-\beta)(1-\epsilon)t\big]_+^{1/(1-\beta)}\end{equation}
We estimate L$g$~in\begin{equation*}
M_1=\big\{(x,t):x_{\epsilon}<x<\eta_{\ell}(t),~~~0<t<\delta_1\big\},\end{equation*}
\begin{equation*}
\eta_{\ell}(t)=-\ell t^{1/(\alpha(1-\beta)}, ~~\ell(\epsilon)=(C+\epsilon)^{-1/\alpha}\big[b(1-\beta)(1-\epsilon)\big]^{1/\alpha(1-\beta)},\end{equation*}
where ~$\delta_1>0$~ is chosen such that~$\eta_{\ell(\epsilon)}(\delta_1)=x_{\epsilon}.$~We have
\begin{equation*}Lg_{\epsilon}=bg_{\epsilon}^{\beta}\{\epsilon+S\}\end{equation*}
\[S=-b^{-1}(p-1)\alpha^{p-1}(\alpha(1-\beta)-1)(C+\epsilon)^{p-1-\beta}(-x)_+^{\alpha(p-1-\beta)-p}\Big\{g_{\epsilon}|x|^{-\alpha}\Big /(C+\epsilon)\Big\}^{\beta( p-2)}\]
\begin{equation*}
-b^{-1}\beta(p-1)\alpha^{p}(C+\epsilon)^{p-1-\beta}(-x)_+^{\alpha(p-1-\beta)-p}\Big\{g_{\epsilon}|x|^{-\alpha}\Big /(C+\epsilon)\Big\}^{\beta( p-1)-1}\end{equation*}
\begin{equation*}
=-b^{-1}\alpha^{p-1}(C+\epsilon)^{p-1-\beta}(-x)_+^{\alpha(p-1-\beta)-p}\Big\{g_{\epsilon}|x|^{-\alpha}\Big/(C+\epsilon)\Big\}^{\beta( p-1)-1}S_1,\end{equation*}
\begin{equation*}
S_1=\Big\{(\alpha(1-\beta)-1)(p-1)\Big[g_{\epsilon}|x|^{-\alpha}\Big/ (C+\epsilon)\Big]^{1-\beta}+\alpha \beta(p-1)\Big\}.\end{equation*}
If~$ \beta(p-1) \geq 1$,~then we can choose~$x_{\epsilon}<0$~such that (with sufficiently small $|x_{\epsilon}|$)
\[|S|<\frac{\epsilon}{2}~~\text{in}~~M_1\]
Thus we have
\[\text{L}g_{\epsilon}~~~~~~~>~~~b(\epsilon/2)g^{\beta}_{\epsilon}\quad (\text{respectively}, \quad \text{L}g_{-\epsilon}<-b\big(\epsilon/2)g^{\beta}_{-\epsilon}\big)~~\text{in}~~M_1,\]
\begin{equation*}
\text{L}g_{\pm \epsilon}~~~~~~=~~~0\quad ~\text{for}~~x>\eta_{\ell(\pm \epsilon)}(t),~~0<t\leq \delta_1,\end{equation*}
\[g_{\epsilon}(x,0)~~\geq~~ u_0(x)\quad\big(\text{respectively}, ~~g_{-\epsilon}(x,0)\leq u_0(x)\big),~~x\geq x_{\epsilon}.~~~~~~~~~\]
Since~$u$~and~$g$~are continuous functions,~$\delta=\delta(\epsilon)\in (0,\delta_1]$~may be chosen such that
\[g_{\epsilon}(x_{\epsilon},t)\geq u(x_{\epsilon},t)~~\big(\text{respectively},~~~ g_{-\epsilon}(x_{\epsilon},t)\leq u(x_{\epsilon},t)\big),~~0\leq t\leq \delta.~~\]
From comparison Lemma~\ref{lemma2} it follows that\begin{subequations}
\begin{equation}\label{CP16a'''}g_{-\epsilon}\leq u\leq g_{\epsilon}\quad x\geq x_{\epsilon},~~0\leq t\leq \delta,\end{equation}
\begin{equation}\label{CP16b'''}\eta_{\ell(-\epsilon)}(t)\leq \eta(t)\leq \eta_{\ell(\epsilon)},~~ 0\leq t\leq \delta, \end{equation}\end{subequations}
which imply  \eqref{CP17'} and  \eqref{CP18'}.

Let~$\beta(p-1)<1$.In this case the left-hand side of \eqref{CP16a'''}, \eqref{CP16b'''} may be proved similarly. Moreover, we can replace~$1+\epsilon$~with $1$ in~$g_{-\epsilon}$~and~$\eta_{\ell}(-\epsilon).$

To prove a relevant upper estimation, consider a function
\[g(x,t)=C_6\big(-\zeta_5t^{\frac{1}{\alpha(1-\beta)}}-x\big)_+^{\alpha}~\text{in}~G_{\ell,\delta},\]
\[G_{\ell,\delta}=\{(x,t):\eta_{\ell}(t)<x<+\infty,~~0<t<\delta\},\]
where~$\ell \in (\ell_*,+\infty)$~and
\begin{equation*}
\zeta_5= (\ell_*/\ell)^{\alpha(1-\beta)}(1-\epsilon)\ell, \ \ 
%\end{equation*}\begin{equation*}
C_6=\big[1-(\ell_*/\ell)^{\alpha(1-\beta)}(1-\epsilon)\big]^{-\alpha}\big[C^{1-\beta}-\ell^{-\alpha(1-\beta)}b(1-\beta)(1-\epsilon))\big]^{1/(1-\beta)},\end{equation*}
From  \eqref{CP18'} it follows that for $\forall \ell>\ell_*$~and for $\forall \epsilon>0$~there exists a~$\delta=\delta(\epsilon,\ell)>0$~such that
\begin{equation}\label{CP17'''}u(\eta_{\ell}(t),t)\leq [C^{1-\beta}\ell^{\alpha(1-\beta)}-b(1-\beta)(1-\epsilon)]^{\frac{1}{1-\beta}}t^{\frac{1}{1-\beta}},~~0\leq t\leq \delta.\end{equation}
Calculating ~L$g$~in
\[G^+_{\ell,\delta}=\{(x,t):\eta_{\ell}(t)<x<-\zeta_5t^{\frac{1}{\alpha(1-\beta)}},~~0<t<\delta\},\]
we have
\[\text{L}g=bg^{\beta}S,~~S=1-(b(1-\beta))^{-1}\zeta_5 C_6^{1/\alpha}\{gt^{1/(\beta-1)}\}^{1-\beta-1/\alpha}
-b^{-1}(\alpha -1)(p-1)\alpha^{p-1}C_6^{p/\alpha}g^{p-1-\beta-(p/\alpha)}.\]
Since
\begin{equation*}S_x\geq 0~~\text{in}~~G^+_{l,\delta},\end{equation*}
\begin{equation*}S\geq S|_{x=\eta_{\ell}(t)}=1-(b(1-\beta))^{-1}\zeta_5 C_6^{1-\beta}(\ell-\zeta_5)^{\alpha(1-\beta)-1}\end{equation*}\begin{equation*}
-b^{-1}(\alpha -1)(p-1)\alpha^{p-1}C_6^{p-1-\beta}\{(\ell-\zeta_5)t^{1/\alpha(1-\beta)}\}^{\alpha(p-1-\beta)-p}.\end{equation*}
Then we have\begin{equation*}
S\geq\epsilon-b^{-1}C_6^{p-1-\beta}(\alpha -1)(p-1)\alpha^{p-1}\{(\ell-\zeta_5)t^{1/\alpha(1-\beta)}\}^{\alpha(p-1-\beta)-p}~\text{in}~G^+_{\ell,\delta}.\end{equation*}
Hence, we can choose~$\delta=\delta(\epsilon)>0$~so small that\begin{subequations}
\begin{equation}\label{CP18a'''}\text{L}g~\geq ~b(\epsilon/2)g^{\beta}~~\text{in}~~G_{\ell,\delta}^+. \end{equation}
Using \eqref{CP17'''}, we can apply Lemma~\ref{lemma1}  in~$G'_{\ell,\delta}=G_{\ell,\delta}\cap \{x<x_0\}$, for $\forall x_0>0$. We have
\begin{equation}\label{CP18b'''}\text{L}g=0~~\text{in}~~G'_{\ell,\delta} \backslash \bar G^+_{\ell,\delta},\end{equation}
\begin{equation}\label{CP18c'''} u(\eta_{\ell}(t),t)\leq \big[C^{1-\beta}\ell^{\alpha(1-\beta)}-b(1-\beta)(1-\epsilon)\big]^{\frac{1}{1-\beta}}t^{\frac{1}{1-\beta}}
=C_6(\ell-\zeta_5)^{\alpha}t^{\frac{1}{(1-\beta)}}=g(\eta_{\ell}(t),t),\;0\leq t\leq \delta.\end{equation}\begin{equation}\label{CP18d'''}u(x_0,t)=g(x_0,t)=0, \ 0\leq t\leq \delta, \ \ 
u(x,0)~=g(x,0)~=0, \ 0\leq x\leq x_0.\end{equation}
\end{subequations}
Since $x_0>0$~is arbitrary, from \eqref{CP18a'''}-\eqref{CP18d'''} and comparison principle it follows that for all ~$\ell>\ell_*$~and $\epsilon>0$~ there exists~$\delta=\delta(\epsilon,\ell)>0$~such that
\begin{equation}\label{CP19'''}u(x,t)\leq C_6(-\zeta_5 t^{\frac{1}{\alpha(1-\beta)}}-x)_+^{\alpha}~~\text{in}~~\bar G_{\ell,\delta.}\end{equation}
Since \eqref{CP18'} is valid along~$x=\eta_{\ell}(t), ~\delta$~may be chosen so small that
\begin{equation}\label{CP20'''}-\ell t^{1/\alpha(1-\beta)}\leq \eta(t)\leq -\zeta_5t^{1/\alpha(1-\beta)},\;0\leq t\leq \delta.\end{equation}
Since~$\ell>\ell_*$~and~$\epsilon>0$~are arbitrary numbers, \eqref{CP17'} follows from \eqref{CP20'''}.

(4a)~~This case is immediate.

(4b)~~Let~$\beta=1, ~\alpha>p/(p-2)~$. As before, from \eqref{CP3}, \eqref{CP5''} follows.  Then consider a function
\[g(x,t)=(C-\epsilon)(-x)_+^{\alpha}\text{exp}(-bt),\]
which satisfies
\[\text{L}g\leq 0~~\text{for}~x_{\epsilon}<x<0,\; t>0;\quad Lg=0~~\text{for}~x>0,\; t>0.\]
We can choose~$\delta=\delta(\epsilon)>0$~such that
\[g(x_{\epsilon},t)\leq u(x_{\epsilon},t),~~~~0\leq t\leq \delta_{\epsilon},\]
and from a comparison principle, the left-hand side of \eqref{CP20'} follows.
To prove the right-hand side, consider
\[g(x,t)=(C+\epsilon)(-x)_+^{\alpha}\text{exp}(-bt)\big[1-\epsilon (b(p-2))^{-1}\big(1-\text{exp}(-b(p-2)t)\big)\big]^{1/2-p}.\]
We have
\begin{equation*}
\text{L}g=(p-2)^{-1}(C+\epsilon)(-x)_+^{\alpha}\text{exp}(-b(p-1)t)g^{p-1}\end{equation*}
\[\times \big\{\epsilon-(p-2)\alpha^{p-1}(\alpha-1)(p-1)(C+\epsilon)^{p-2}(-x)_+^{\alpha(p-2)-p}\big\},~~x<0,\; t>0,\]
and hence, if $|x_{\epsilon}|$ is small enough,
\[\text{L}g\geq 0 ~~~\text{for}~~~x_{\epsilon}<x<0,\; t>0; \quad \text{L}g=0~~~ \text{for}~~~x>0,~t>0.\]
As before, a comparison principle implies the right-hand side of \eqref{CP20'}. The estimations \eqref{CP21'}-\eqref{CP23'} in the cases (4c) and (4d) may be proved similarly.

(II)~$b=0$.

(1)~Let~$p>2, ~0<\alpha<p/(p-2)$.\\
First assume that~$u_0$~is defined by \eqref{CP4}. The self-similar form \eqref{CP25'} and the formula\eqref{CP26'} are well-known results (see Lemma~\ref{lemma2}). To prove \eqref{CP27'}, consider a function
\[g(x,t)=t^{\alpha/(p-\alpha(p-2))}f(\xi).\]
We have
\[\text{L}g=t^{(\alpha(p-1)-p)/(p-\alpha(p-2))}\mathcal{L}_tf,\]
\[\mathcal{L}_tf=\frac{\alpha}{p-\alpha(p-2)}f-\frac{1}{p-\alpha(p-2)}\xi f'-\big(|f'|^{p-2}f'\big)'.\]
Choose
\[f(\xi)=C_0(\xi_0-\xi)_+^{(p-1)/(p-2)},~~~~0<\xi<+\infty\]
where $C_0$ and $\xi_0$ are some positive constants. Then we have
\[\mathcal{L}_tf=(p-\alpha(p-2))^{-1}(p-1)(p-2)^{-1}C_0(\xi_0-\xi)^{1/p-2}R(\xi) \quad \text{for}~ 0\leq \xi \leq \xi_0,~~~~t>0\]
\[R(\xi)=\alpha(p-2)(p-1)^{-1}\xi_0+(1-\alpha(p-2)(p-1)^{-1})\xi-(p-1)^{p-1}(p-2)^{-(p-1)}(p-\alpha(p-2))C_0^{p-2}\]
To prove an upper estimation we take ~$C_0=C_5,~~\xi_0=\xi_4$. Then we have
\[R(\xi)\geq\nu_{\alpha}\xi_4-(p-1)^{p-1}(p-2)^{-(p-1)}(p-\alpha(p-2))C_5^{p-2}=0\quad \text{for}~~0\leq \xi \leq \xi_4,\]
where
\[\nu_{\alpha}=\{1~~\text{if}~~\alpha \geq (p-1)(p-2)^{-1};~\alpha(p-2)(p-1)^{-1}~~\text{if}~\alpha<(p-1)(p-2)^{-1}\}.\]
Hence
\[\text{L}g~\geq 0~~\text{for}~~0<x<\xi_4 t^{1/p-\alpha(p-2)},~~~t>0,\]
\[\text{L}g~= 0~~\text{for}~~0>\xi_4 t^{1/p-\alpha(p-2)},~~~~~~~~~t>0,\]
\[u(0,t)=g(0,t),~t\geq 0;~u(x,0)=g(x,0),~x\geq 0\]
and a comparison principle imply the right-hand side of \eqref{CP27'}. The left-hand side of \eqref{CP27'} may be established similarly if we take $C_0=C_4,~~\xi_0=\xi_3$. \eqref{CP2'} and  \eqref{24ab} follow from Lemma~\ref{lemma2}. Finally, \eqref{CP5'}-\eqref{CP7'} easily follow from \eqref{CP26'} and \eqref{CP27'}. If $u_0$ satisfies \eqref{CP3} with $0<\alpha<p/(p-2),$ then \eqref{CP1'}-\eqref{CP3'} follow from Lemma~\ref{lemma2}.

The cases (2) and (3) are immediate.

%%%%%%%%%%%%%%%%%%%%%%%%%%%%%%%%%%

\newpage

\textbf{Appendix.} We give here explicit values of the constants used in  section 2 in the outline of the results for Case (I(2)) and later in section 4 during the proof of these results.\\\\
$\zeta_1=A^{\frac{p-2}{p}}(1-\beta)^{\frac{1}{p}}(p-1)\big(1+b(1-\beta)A_1^{\beta-1}\big)^{-\frac{1}{p}}(p-2)^{-1},~~~~$\\

$C_1=A_1 \zeta_1^{-\mu}~~~~~~~~\text{if}~~~~\beta(p-1)>1$,\\\\
$\zeta_1=A_1^{\frac{p-2}{p}}\big((1-\beta)(1+\beta)p^{p-1}(p-1)\big)^{\frac{1}{p}}\big(1+b(1-\beta)A_1^{\beta-1}\big)^{-\frac{1}{p}}(p-1-\beta)^{-1},~~~$\\

$C_1=A_1 \zeta_1^{-\frac{p}{p-1-\beta}},~~~~~~~~\text{if}~~~~\beta(p-1)<1$,\\\\
$\zeta_2=A_1^{\frac{p-2}{p}}\big((1-\beta)(1+\beta)p^{p-1}(p-1)\big)^{\frac{1}{p}}\big(1+b(1-\beta)A_1^{\beta-1}\big)^{-\frac{1}{p}}(p-1-\beta)^{-1},~~~~$\\

$ C_2=A_1 \zeta_2^{-\frac{p}{p-1-\beta}},~~~~~\text{if}~~~~~~~\beta(p-1)>1$,\\\\
$\zeta_2=\Big(A_1/C_*\Big)^{\frac{p-1-\beta}{p}},~~~~~C_2=C_* ,~~~~~~~\text{if}~~~\beta(p-1)<1$,\\\\
$\bar \zeta_2= A_1^{\frac{p-2}{p}} \Big(\frac{p(p-1)^{p}(p-2)^{1-p}(1-\beta)}{p(p-2)-\beta(p-1)+1}\Big)^{\frac{1}{p}},~~~~~~~~\bar C_2=A_1 \bar \zeta_2^{-\frac{(p-1)}{(p-2)}},~~~~~~~\text{if}~~~~\beta(p-1)>1,$\\\\
$\bar\zeta_2=A^{\frac{p-2}{p}}(1-\beta)^{\frac{1}{p}}(p-1)\big(1+b(1-\beta)A_1^{\beta-1}\big)^{-\frac{1}{p}}(p-2)^{-1},$\\\\
$~~~~~~~~\bar C_2=A_1 \bar \zeta_2^{-\frac{(p-1)}{(p-2)}},~~~~~~~\text{if}~~~~\beta(p-1)<1,$\\\\
$\ell_0=C_*^{\frac{1+\beta-p}{p}}(C_*/C\Big)^{\frac{(1-\beta)(p-1-\beta)}{1-\beta(p-1)}}(b(1-\beta)\theta_*)^{\frac{p-1-\beta}{p(1-\beta)}}$,\\\\
$\zeta_3=C_*^{\frac{1+\beta-p}{p}}\Big[(C_*/C\Big)^{\frac{(1-\beta)(p-1-\beta)}{1-\beta(p-1)}}-1\Big](b(1-\beta)\theta_*)^{\frac{p-1-\beta}{p(1-\beta)}}$,\\\\
$\theta_*=\Big[1-\Big(C/C_*\Big)^{p-1-\beta}\Big]\Big[\Big(C_*/C\Big)^{\frac{(1-\beta)(p-1-\beta)}{1-\beta(p-1)}}-1\Big]^{-1},$\\\\
$\ell_1=C^{\frac{1+\beta-p}{p}}\Big[b(1-\beta)(\delta_* \Gamma)^{-1}\Big((1-\delta_* \Gamma)-\big(1-\delta_* \Gamma\big)^{1-p}\Big(C/C_*\Big)^{p-1-\beta}\Big) \Big]^{\frac{p-1-\beta}{p(1-\beta)}}$,\\\\
$\zeta_4=\delta_* \Gamma \ell_1,~~~~~~~~\Gamma=1- (C/C_*)^{\frac{p-1-\beta}{p}},~~~~~~~~C_3=C \big(1-\delta_* \Gamma\big)^{\frac{p}{1+\beta-p}},$~\\\\\
where ~$\delta_*\in(0,1)$~satisfies

$g(\delta_*)=\underset{[0;1]}\max g(\delta),\quad ~~~g(\delta)=\delta^{\frac{1+\beta(1-p)}{p(1-\beta)}}\Big[(1-\delta \Gamma)-\Big(C/C_*\Big)^{p-1-\beta}\big(1-\delta \Gamma\big)^{1-p}\Big) \Big],$\\
$\ell_*=C^{-\frac{1}{\alpha}}\big(b(1-\beta)\big)^{1/(\alpha(1-\beta))}$\\\\
$\zeta_5= (\frac{\ell_*}{\ell})^{\alpha(1-\beta)}(1-\epsilon)\ell,~~~~~~~~~\text{if}~~\beta(p-1)<1$,\\\\
$C_6=\big(1-(\frac{\ell_*}{\ell})^{\alpha(1-\beta)}(1-\epsilon)\big)^{-\alpha}\big[C^{1-\beta}-\ell^{-\alpha(1-\beta)}b(1-\beta)(1-\epsilon))\big]^{\frac{1}{1-\beta}}.$\\\\

\end{document}